\documentclass[11pt]{amsart}

\usepackage{amsmath, amsthm}
\usepackage{amssymb}
\usepackage{amsfonts}
\usepackage{latexsym}
\usepackage{mathrsfs}
\usepackage{bm}
\usepackage{graphicx}
\usepackage[dvipsnames,usenames]{color}

\newcommand{\be}{\begin{equation}}
\newcommand{\ene}{\end{equation}}
\newcommand{\br}{\begin{remark}}
\newcommand{\er}{\end{remark}}
\newcommand{\bl}{\begin{lemma}}
\newcommand{\el}{\end{lemma}}
\newcommand{\bcor}{\begin{cor}}
\newcommand{\ecor}{\end{cor}}
\newcommand{\bpro}{\begin{proposition}}
\newcommand{\epro}{\end{proposition}}
\newcommand{\ben}{\begin{enumerate}}
\newcommand{\een}{\end{enumerate}}
\newcommand{\bp}{\begin{proof}}
\newcommand{\ep}{\end{proof}}
\newcommand{\bpo}{\begin{pro}}
\newcommand{\epo}{\end{pro}}
\newcommand{\beq}{\begin{equation*}}
\newcommand{\eeq}{\end{equation*}}
\newcommand{\bear}{\begin{eqnarray}}
\newcommand{\eear}{\end{eqnarray}}
\newcommand{\beqar}{\begin{eqnarray*}}
\newcommand{\eeqar}{\end{eqnarray*}}
\newcommand{\brem}{\begin{remark*}}
\newcommand{\erem}{\end{remark*}}
\newcommand{\bt}{\begin{theorem}}
\newcommand{\et}{\end{theorem}}
\newcommand{\bDe}{\begin{Def}}
\newcommand{\eDe}{\end{Def}}

\newcommand{\abs}[1]{\lvert #1 \rvert}

\newcommand{\Z}{\mathbb{Z}}
\newcommand{\R}{\mathbb{R}}

\newcommand{\bH}{\mathbb{H}}
\newcommand{\Cc}{\mathcal{C}}
\newcommand{\N}{\mathcal{N}}

\newcommand{\grad}{\nabla}

\DeclareMathOperator{\im}{Im}
\DeclareMathOperator{\re}{Re}
\renewcommand{\Re}{\re}
\renewcommand{\Im}{\im}
\DeclareMathOperator{\sys}{sys}

\renewcommand{\tilde}[1]{\widetilde{#1}}

\newcommand{\eps}{\varepsilon}

\newcommand{\wep}{Weil-Petersson }

\newcommand{\sT}{\mathcal{T}}

\newcommand{\sM}{\mathcal{M}}
\newcommand{\lsys}{\ell_{sys}(X)}

\renewcommand{\leq}{\leqslant}
\renewcommand{\geq}{\geqslant}

\newcommand{\RS}{Riemann surface}

\newcommand{\WP}{Weil-Petersson}

\newcommand{\nbhd}{neighbourhood}

\newcommand{\cts}{continuous}

\DeclareMathOperator{\MC}{\mathcal{MC}}

\DeclareMathOperator{\Area}{Area}
\DeclareMathOperator{\Vol}{Vol_{WP}}
\DeclareMathOperator{\diam}{diam}
\DeclareMathOperator{\Inrad}{InRad}

\DeclareMathOperator{\dist}{dist}
\DeclareMathOperator{\inj}{inj}

\DeclareMathOperator{\arcsinh}{arcsinh}

\DeclareMathOperator{\Teich}{Teich}

\DeclareMathOperator{\Diff}{Diff}

\newcommand{\Mod}{\mbox{\rm Mod}}

\DeclareMathOperator{\dVol}{dVol}
\DeclareMathOperator{\dArea}{dArea}

\newcommand{\p}{\partial}

\newcommand{\param}{{\mathchoice{\mkern1mu\mbox{\raise2.2pt\hbox{$\centerdot$}}\mkern1mu}{\mkern1mu\mbox{\raise2.2pt\hbox{$\centerdot$}}\mkern1mu}{\mkern1.5mu\centerdot\mkern1.5mu}{\mkern1.5mu\centerdot\mkern1.5mu}}}

\numberwithin{equation}{section}

\theoremstyle{plain}
\newtheorem{theorem}{Theorem}[section]

\newtheorem{corollary}[theorem]{Corollary}
\newtheorem{lemma}[theorem]{Lemma}
\newtheorem{proposition}[theorem]{Proposition}

\newtheorem*{rem*}{Remark}
\newtheorem*{def*}{Definition}

\theoremstyle{definition}

\newtheorem{question}{Question}

\newtheorem{remark}[theorem]{Remark}
\theoremstyle{definition}
\newtheorem*{remarksenv}{Remarks}
\newtheorem*{rem}{Remark}

\begin{document}


\date{\today}

\title[Uniform Lower Bound]
{A new uniform lower bound on Weil-Petersson distance}

\author{Yunhui Wu}
\address{Tsinghua University, Haidian District, Beijing 100084, China}
\email{yunhui\_wu@mail.tsinghua.edu.cn}


\maketitle

\begin{abstract}
In this paper we study the injectivity radius based at a fixed point along Weil-Petersson geodesics. We show that the square root of the injectivity radius based at a fixed point is $ 0.3884$-Lipschitz on Teichm\"uller space endowed with the Weil-Petersson metric. As an application we reprove that the square root of the systole function is uniformly Lipschitz on Teichm\"uller space endowed with the Weil-Petersson metric, where the Lipschitz constant can be chosen to be $0.5492$. Applications to asymptotic geometry of moduli space of Riemann surfaces for large genus will also be discussed.   
\end{abstract}

\section{Introduction}
Let $S_g$ be a closed surface of genus $g$ $(g\geq 2)$, and $\sT_g$ be the Teichm\"uller space of $S_g$. Let $\Teich(S_g)$ be the space $\sT_g$ endowed with the \WP \ metric. The mapping class group $\Mod(S_g)$ acts on $\Teich(S_g)$ by isometries. The moduli space $\sM_g$ of $S_g$ endowed with the \WP \ metric, is realized as the quotient $\Teich(S_g)/\Mod(S_g)$. Let $\sM_{-1}$ be the space of complete Riemannian metrics on $S_g$ of constant Gauss curvature $-1$. It is known that $\sT_g=\sM_{-1}/\Diff_0(S_g)$ where $\Diff_0(S_g)$ is the group of diffeomorphisms of $S_g$ isotopic to the identity. Let $p \in S_g$ be fixed and $\tilde{X}\in \sM_{-1}$ be a hyperbolic metric on $S_g$. The \emph{injectivity radius} $\inj_{\tilde{X}}(p)$ of $\tilde{X}$ at $p$ is half of the length of a shortest nontrivial geodesic loop based at $p$. The geodesic loop based at $p$ realizing $\inj_{\tilde{X}}(p)$ may not be unique. It is known that $\inj_{\tilde{X}}(p)$ is bounded from above by a positive constant only depending on $g$. Let 
\[\pi:\sM_{-1}\to \sT_{g}\]
be the natural projection. In order to transfer the quantity $\inj_{(\cdot)}(p)$ onto \WP \ geodesics in $\sT_g$, we make the following definition. First we recall in \cite{RT18,RT18-JDG} that for a smooth path $c(t) \subset \sM_{-1}$, we say $c(t)$ is a \emph{horizontal curve} if for each $t$, there exists a holomorphic quadratic differential $\Psi(t)$ of $c(t)$ such that the variation of hyperbolic metrics satisfies $\frac{\partial c(t)}{\partial t}=\Re{\Psi(t)}$. A smooth path in $\sT_g$ can always be lifted onto a horizontal curve in $\sM_{-1}$. Throughout this paper we always assume that parameters are proportional to arc-length parameters for both geodesics in hyperbolic surfaces and smooth \WP \ paths in Teichm\"uller space of \RS s. It is not well-defined for $\inj_{(\cdot)}(p)$ on $\sT_g$ because point $p$ is clearly not invariant by diffeomorphisms of $S_g$. To fix this problem, we first pick the \WP \ geodesic joining $X$ and $Y$, and lift this \WP \ geodesic onto a horizontal curve $c:[0,1]\to \sM_{-1}$. Then we consider the injectivity radius function $\inj_{c(t)}(p)$ along $c([0,1])$ which is well-defined. Now we define
\begin{def*}
Fix a point $p \in S_g$ $(g\geq 2)$. For any $X, Y\in \sT_g$, we define
\[\left|\sqrt{\inj_X(p)}-\sqrt{\inj_Y(p)}\right|:=\sup_{c}\left|\sqrt{\inj_{c(0)}(p)}-\sqrt{\inj_{c(1)}(p)}\right|\]
where $c:[0,1] \to \sM_{-1}$ runs over all smooth horizontal curves with $\pi(c(0))=X$, $\pi(c(1))=Y$ and $\pi(c([0,1]))\subset \sT_g$ is the \WP \ geodesic joining $X$ and $Y$. 
\end{def*}  

\noindent The definition above actually does not depend on the choice of $p$. One may see the following remark for an equivalent definition.

\begin{rem*}
It is known that any two horizontal lifts in $\sM_{-1}$ of a smooth curve in $\sT_g$ differ by an element in $\Diff_0(S_g)$ (see \emph{e.g.} \cite[Chapter $2$]{Tromba-book}), and the group $\Diff_0(S_g)$ acts transitively on $S_g$. So the definition above is equivalent to 
\[\left|\sqrt{\inj_X(p)}-\sqrt{\inj_Y(p)}\right|:=\sup \limits_{q\in S_g}\left|\sqrt{\inj_{c'(0)}(q)}-\sqrt{\inj_{c'(1)}(q)}\right|\]
where $c':[0,1] \to \sM_{-1}$ is a horizontal lift of the \WP \ geodesic joining $X$ and $Y$. 
\end{rem*}

In this paper, we show that
\bt \label{mt}
Fix a point $p \in S_g$ $(g\geq2)$. Then for any $X, Y \in \sT_g$, 
$$\left|\sqrt{\inj_X(p)}-\sqrt{\inj_Y(p)}\right|\leq  0.3884 \dist_{wp}(X,Y)$$
where $\dist_{wp}$ is the \WP \ distance.
\et

\begin{rem}
Rupflin and Topping in \cite[Section 2]{RT18} showed that $$\left|\sqrt{\inj_X(p)}-\sqrt{\inj_Y(p)}\right|\leq c(g) \dist_{wp}(X,Y)$$ where $c(g)>0$ is a constant depending on $g$. Our approach is similar to that of Rupflin and Topping, but using a detailed analysis of injectivity radius along shortest geodesic loops and a recent uniform bound for harmonic Beltrami differentials on thin parts \cite{BW-curv}, we are able to obtain the above uniform bound independent of $g$. The Lipschitz constant $0.3884$ above is not optimal. More refined arguments in this proof can improve this uniform constant. In general, it is difficult to measure the \WP \ distance on $\sT_g$. One may see \cite{Brock03, BBB19, BB-20, CP12, KM18, RT18, Schl13, Schl19, Wolpert08, W-inradius} for related bounds on \WP \ distances.
\end{rem}

Recall that the \emph{systole} $\ell_{sys}(X)$ of $X\in \sT_g$ is the length of a shortest nontrivial closed geodesic in $X$. Which is also the same as $2\min_{p\in \tilde{X}}\inj_{\tilde{X}}(p)$ where $\tilde{X}\in \sM_{-1}$ is a hyperbolic metric on $S_g$ with $\pi(\tilde{X})=X\in \sT_g$. As a direct application of Theorem \ref{mt}, we prove
\begin{corollary} \label{mt-2}
For any $X, Y \in \sT_g$ $(g\geq 2)$, 
$$\left|\sqrt{\ell_{sys}(X)}-\sqrt{\ell_{sys}(Y)}\right|\leq 0.5492 \dist_{wp}(X,Y).$$
\end{corollary}
\bp
Without loss of generality, one may assume that $$\ell_{sys}(X)\geq \ell_{sys}(Y).$$ Let $c:[0,1] \to \sM_{-1}$ be a horizontal curve with $\pi(c(0))=X$, $\pi(c(1))=Y$ and $\pi(c([0,1]))\subset \sT_g$ is the \WP \ geodesic joining $X$ and $Y$. Let $\alpha \subset c(1)$ be a shortest closed geodesic. So for any $p \in \alpha$, we have
\[2\inj_{c(1)}(p)=\ell_{sys}(Y) \quad \emph{and} \quad  2\inj_{c(0)}(p)\geq \ell_{sys}(X).\] 
Then by Theorem \ref{mt} we get
\[\sqrt{\ell_{sys}(X)}-\sqrt{\ell_{sys}(Y)}\leq \sqrt{2\inj_{c(0)}(p)}-\sqrt{2\inj_{c(1)}(p)}\leq 0.5492\dist_{wp}(X,Y)\]
as desired.
\ep

\begin{rem}
\ben
\item It was shown in \cite{W-inradius} that \[\left|\sqrt{\ell_{sys}(X)}-\sqrt{\ell_{sys}(Y)}\right|\leq K \dist_{wp}(X,Y)\] where $K>0$ is a uniform (implicit) constant independent of $g$. 
\item Very recently, Bridgeman-Bromberg in \cite{BB-20} show that the uniform constant $K$ above can be chosen to be $\frac{1}{2}$ by a completely different method. 
\een
\noindent Both the proofs in \cite{W-inradius} and \cite{BB-20} rely on certain uniform bound for the \WP \ norm $||\grad(\ell_{\alpha}(X))||_{wp}$ of the \WP \ gradient $\grad(\ell_{\alpha}(X))$ of the geodesic length function $\ell_{\alpha}(\cdot)$ on Teichm\"uller space, where $\alpha \subset X$ is a systolic curve (one may also see \cite{W-sys} for a different proof). In this paper, our proof is totally different without any estimation on $||\grad(\ell_{\alpha}(X))||_{wp}$. Moreover, we are able to obtain the explicit Lipschitz constant above to be $0.5492$, which can be improved by more careful arguments for the proof of Theorem \ref{mt}.  
\end{rem}

The \WP \ completion of the moduli space $\sM_g$ is compact which is homeomorphic to the Deligne-Mumford compactification of the moduli space of \RS s. In particular, the moduli space $\sM_g$ has finite \WP \ diameter and inradius. Cavendish-Parlier \cite{CP12} showed that for large genus the ratio $\frac{\diam(\sM_g)}{\sqrt{g}}$ is bounded below by a uniform positive constant and above by a uniform constant multiple of $\ln(g)$. It is an \emph{open} problem that whether the \wep diameter $\diam(\sM_g)$ of $\sM_g$ is uniformly comparable to $\sqrt{g}\ln(g)$. Recall that the \WP \ \emph{inradius} $\Inrad(\sM_g)$ of $\sM_g$ is defined as
\[\Inrad(\sM_g):=\max_{X\in \sM_g}\dist_{wp}(X, \p \sM_g)\]
where $\p \sM_g$ is the boundary of $\sM_g$ consisting of nodal surfaces. It was shown in \cite{W-inradius} that as $g\to \infty$, the \WP \ inradius $\Inrad(\sM_g)$ is uniformly comparable to $\sqrt{\ln(g)}$. More precisely, there exists a uniform constant $K'>0$ independent of $g$ such that $$K'\leq \frac{\Inrad(\sM_g)}{\sqrt{\ln(g)}} \leq \sqrt{4\pi}$$ where the uniform (implicit) constant $K'$ depends on the work of Buser-Sarnak \cite{BS94}. The following question is natural and interesting.
\begin{question}\label{ques-1}
Does $\lim\limits_{g\to \infty}\frac{\Inrad(\sM_g)}{\sqrt{\ln(g)}}$ exist? If exists, what is its value?  
\end{question} 

Set
\[\sys(g)=\max_{X\in \sM_g}\lsys.\]
It is known that $$\lsys\leq 2\ln(4g-2)$$ for all $X\in \sM_g$. Buser and Sarnak in \cite{BS94} showed that $$\sys(g)\geq U\ln(g)$$ for some uniform constant $U>0$. Moreover, they also showed that there exists a sequence $\{g_k\}_{k\geq 1}$ of positive integers tending to infinity such that for each $g_k$, there exists a closed hyperbolic surface $\mathcal{X}_{g_k}$ of genus $g_k$ with
\[\ell_{sys}(\mathcal{X}_{g_k})\geq \frac{4\ln(g_k)}{3}-U'\]
where $U'>0$ is a uniform constant independent of $g$. Thus, the quantity $\sys(g)$ is uniformly comparable to $\ln(g)$ as $g\to \infty$. Moreover 
$$\frac{4}{3}\leq \limsup_{g\to \infty}\frac{\sys(g)}{\ln(g)}\leq 2.$$ 
By applying the proof of Theorem \ref{mt} it is not hard to see that
$$2.0472\leq \liminf_{g\to \infty}\frac{\Inrad(\sM_g)}{\sqrt{\sys(g)}}\leq 2.5066$$
and
$$2.3696\leq \limsup_{g\to \infty}\frac{\Inrad(\sM_g)}{\sqrt{\ln(g)}}\leq 3.5449.$$

In this paper we show that
\bt\label{mt-3}
The following limit holds:
\[\lim_{g\to \infty}\frac{\Inrad(\sM_g)}{\sqrt{\sys(g)}}=\sqrt{2\pi}\sim 2.5066.\]
\et

\begin{rem*}
\ben
\item Theorem \ref{mt-3} was firstly obtained in \cite{BB-20} by Bridgeman and Bromberg. We are grateful to M. Bridgeman for kindly sharing their latest version of \cite{BB-20}.

\item The proof of \cite{BB-20} relies on bounds of $||\grad(\ell_\alpha(X))||_{wp}$ in terms of functions on collars. Our proof is a slightly refined argument of the proof of \cite[Theorem  1.1]{W-inradius} where we bound $||\grad(\ell_\alpha(X))||_{wp}$ in terms of functions on hyperbolic disks. In both cases, $\alpha$ is a systolic curve of $X$. 
\een
\end{rem*}
\noindent Theorem \ref{mt-3} reduces Question \ref{ques-1} to study the following one which has no metric involved on $\sM_g$.
\begin{question}\label{ques-2}
Does $\lim\limits_{g\to \infty}\frac{\sys(g)}{\ln(g)}$ exist? If exists, what is its value?  
\end{question}

Recall that for any hyperbolic surface $X\in \sM_g$, the \emph{Bers' constant} $\mathcal{B}_g(X)$ at $X$ is the smallest positive number such that there exist $(3g-3)$ disjoint simple closed geodesics $\{\gamma_i\}_{i=1}^{3g-3}$ on $X$ with $$\max_{1\leq i \leq (3g-3)}\ell_{\gamma_i}(X)\leq \mathcal{B}_g(X).$$ It is known \cite[Chapter 5]{Buser10} that $$\sqrt{6g}-2\leq \sup_{X\in \sM_g}\mathcal{B}_g(X)\leq 26(g-1).$$ Buser \cite{Buser10} \emph{conjectures} that $\sup_{X\in \sM_g}\mathcal{B}_g(X)$ is uniformly comparable to $\sqrt{g}$. Fix any $L>0$, we set the subset $\MC(\leq L)\subseteq \sM_g$ as
\[\mathcal{MC}(\leq L):=\left\{X\in \sM_g; \ \mathcal{B}_g(X)\leq L\right\}.\]
So $\mathcal{MC}\left(\leq \sup_{X\in \sM_g}\mathcal{B}_g(X)\right)=\sM_g$. Let $L_g=\epsilon \ln(g)$ for large $g$ and small enough $\epsilon>0$. By applying \cite{Mirz13} and Theorem \ref{mt}, we show that for large enough $g$, the $1.0511\sqrt{\ln(g)}$-\nbhd \ of $\mathcal{MC}(\leq L_g)$ can be arbitrarily small in $\sM_g$ in the following sense.
\bt \label{mt-4}
For any small enough $\epsilon>0$ and if $L_g=\epsilon \ln(g)$, then 
\[\lim_{g\to \infty}\frac{\Vol\left(\left\{X\in \sM_g; \ \dist_{wp}\left(X,\mathcal{MC}(\leq L_g)\right)> 1.0511\sqrt{\ln(g)}\right\}\right)}{\Vol(\sM_g)}=1\]
where $\Vol(\cdot)$ is the \WP \ volume.
\et

\noindent \textbf{Plan of the paper.} Section \ref{sec-prel} provides some necessary background and the basic properties on two-dimensional hyperbolic geometry and Teichm\"uller theory. In Section \ref{sec-two bounds} we prove two bounds for the injectivity radius along a shortest geodesic loop based at a fixed point. A technical inequality is provided in Section \ref{sec-ine} which is crucial in the proof of Theorem \ref{mt}. In Section \ref{sec-mt} we finish the proof of Theorem \ref{mt}. Theorem \ref{mt-3} is shown in Section \ref{mg-large}. And we prove Theorem \ref{mt-4} in Section \ref{sec-mt4}.\\

\noindent \textbf{Acknowledgements.}
The author would like to thank Martin Bridgeman, Ran Ji and Scott Wolpert for helpful conversations on this paper, and thank to Melanie Rupflin, Peter Topping and Shing-Tung Yau for their interests. He especially would like to thank Peter Topping for useful discussions on the Definition on Page $2$. He is also very grateful to anonymous referees for their helpful comments. This work is supported by the NSFC grant No. $12171263$ and a grant from Tsinghua University.


\section{Preliminaries}\label{sec-prel}
In this section we will set up the notations and provide some necessary background on two-dimensional hyperbolic geometry and Teichm\"uller theory of \RS s.

\subsection{Injectivity radius at a point}
Let $X$ be a closed hyperbolic surface. Since the curvature of $X$ is $-1$, the conjugate radius at any point of $X$ is infinity. Thus for any point $p \in X$, the \emph{injectivity radius} $\inj_X(p)$ of $X$ at $p$ is half of the length of a shortest nontrivial geodesic loop based at $p$. Let $$\sigma:[0, 2\inj_X(p)]\to X$$ be such a shortest geodesic loop with $\sigma(0)=\sigma(2\inj_X(p))=p$ of arc-length parameter. Then 
\ben
\item the restriction $\sigma:[0,\inj_X(p)]\to X$ is a minimizing geodesic;

\item the restriction $\sigma:[\inj_X(p),2\inj_X(p)]\to X$ is also a minimizing geodesic.
\een 
For any $r>0$, we let $$B(p;r):=\left\{q\in X; \ \dist(q,p)<r  \right\}$$
be the open geodesic ball centered at $p$ of radius $r$. The open ball $B(p;\inj_X(p))$ is an embedded hyperbolic open disk of radius $\inj_X(p)$. By the Gauss-Bonnet formula we know that $\Area(X)=4\pi(g-1)$. Thus, $$\Area(B(p;\inj_X(p)))=2\pi \left(\cosh(\inj_X(p))-1\right)\leq 4\pi(g-1)$$ which implies that for any $p \in X$,
\bear\label{up-inj}
\inj_X(p)\leq \ln(4g-2).
\eear
We remark here that for all $g\geq 2$, Buser and Sarnak in \cite{BS94} constructed a closed surface $\mathcal{X}_g$ of genus $g$ such that 
\[\inf_{p \in \mathcal{X}_g}\inj_{\mathcal{X}_g}(p)\geq U\ln(g)\]
for some uniform constant $U>0$ independent of $g$.\\
\vspace{-0.1in}

\subsection{Teichm\"uller space and \wep metric.} \label{sec:wp background}
We denote by $S_{g}$ an oriented closed surface of genus $g$ $(g\geq 2)$. The Uniformization Theorem implies that the surface $S_g$ admits hyperbolic metrics of constant curvature $-1$. We let $\sT_g$ be the Teichm\"uller space of surfaces of genus $g$, which we consider as the equivalence classes under the action of the group $\Diff_0(S_g)$ of diffeomorphisms isotopic to the identity of the space of hyperbolic surfaces $X=(S_g,\sigma(z)|dz|^2)$. The tangent space $T_X\sT_g$ at a point $X=(S_g,\sigma(z)|dz|^2)$ is identified with the space of {\it harmonic Beltrami differentials} on $X$, i.e., forms on $X$ expressible as 
$\mu=\overline{\psi}/\sigma$ where $\psi \in Q(X)$ is a holomorphic quadratic differential on $X$. The pointwise norm $|\mu(\cdot)|:X\to \R^{\geq 0}$ gives a \cts \ nonnegative function on $X$. Let $z=x+iy$ and $\dArea=\sigma(z)dxdy$ be the volume form. The \textit{Weil-Petersson metric} is the Hermitian
metric on $\sT_g$ arising from the the \textit{Petersson scalar  product}
\begin{equation}
 \left<\varphi,\psi \right>= \int_X \frac{\varphi \cdot \overline{\psi}}{\sigma^2}\dArea\nonumber
\end{equation}
via duality. We will concern ourselves primarily with its Riemannian part $g_{WP}$. Throughout this paper we denote by $\Teich(S_g)$ the Teichm\"uller space endowed with the Weil-Petersson metric. By definition it is easy to see that the mapping class group $\Mod(S_g)$ acts on $\Teich(S_g)$ as isometries. Thus, the \wep metric descends to a metric, also called the \wep metric, on the moduli space of Riemann surfaces $\sM_g$ which is defined as $\sT_g/\Mod(S_g)$. Throughout this paper we also denote by $\sM_g$ the moduli space endowed with the Weil-Petersson metric.  One may refer to \cite{IT92,Wolpert-book} for more details on \wep geometry. 

\subsection{Uniform bounds on harmonic Beltrami differentials} In this subsection we recall two uniform bounds on the pointwise norm of any harmonic Beltrami differential in terms of the injectivity radius at a point. We first refer to a function $C(r)$ introduced by Teo in \cite{Teo09} which is given by
\bear \label{defn:C-inj}
C(r)=\left(\frac{4\pi}{3}\left(1-\left(\frac{4e^{r}}{(1+e^{r})^2}\right)^3\right)\right)^{-\frac{1}{2}}.
\eear

\noindent It follows that $C(r)$ is decreasing with respect to $r$ and as $r$ tends to zero we have
$$C(r) = \frac{1}{\sqrt{\pi}r} + O(1) .$$ Furthermore $C(r)$ tends to $\sqrt{\frac{3}{4\pi}}$ as $r$ tends to infinity. 
The following property follows by a Taylor expansion of $\mu$ on a hyperbolic disk of radius $r>0$. 
\bpro[Teo, {\cite[Prop 3.1]{Teo09} or \cite[Prop 2.10]{Wolf-W-1}}] \label{bound-teo}
Let $X$ be a closed hyperbolic surface and $\mu$ be a harmonic Beltrami differential on $X$. Then for any $p \in X$,
 $$|\mu(p)|^2  \leq \left(C(\inj_X(p))\right)^2 \int_{B(p;r)}{|\mu(z)|^2 \dArea(z)}, \ \forall \ 0<r\leq \inj_X(p)$$
where the constant $C(\cdot)$ is given by \eqref{defn:C-inj}.
\epro

Proposition \ref{bound-teo} is useful when the injectivity radius at a point is uniformly bounded from below, especially as the injectivity radius goes to infinity. For the case that the injectivity radius at a point is small, we will use the following recent result, which follows by a detailed analysis on the Fourier expansion of $\mu$ on a collar of a short closed geodesic. More precisely, 
\bpro[Bridgeman-Wu, {\cite[Prop 1.1]{BW-curv}}] \label{bound-BW}
Let $X$ be a closed hyperbolic surface and $\mu$ be a harmonic Beltrami differential on $X$. Then for any $p \in X$ with $\inj_X(p)\leq \arcsinh(1)$,
 $$|\mu(p)|^2  \leq \frac{\int_X |\mu(z)|^2 \dArea(z)}{\inj_X(p)}.$$
\epro

\section{Two bounds on injectivity radius}\label{sec-two bounds}
Let $X$ be a closed hyperbolic surface of genus $g\geq 2$ and $\gamma \subset X$ be a non-trivial simple loop. There always exists a unique closed geodesic, still denoted by $\gamma$, representing this loop. The Collar Lemma says that it has a tubular neighborhood which is a topological cylinder with a standard hyperbolic metric. And the width of this cylinder, only depending on the length of $\gamma$, goes to infinity as the length of $\gamma$ goes to $0$. First we recall the following version of the Collar Lemma which will be applied.

\begin{theorem}\cite[Theorem 4.1.1]{Buser10} \label{collar}
	Let $\gamma_1 , \gamma_2, ..., \gamma_m$ be disjoint simple closed geodesics on a closed hyperbolic Riemann surface $X$ of genus $g$, and $\ell(\gamma_i)$ be the length of $\gamma_i$. Then $m\leq 3g-3$ and we can define the collar of $\gamma_i$ by
	$$\Cc(\gamma_i)=\{x\in X_g; \ \dist(x,\gamma_i)\leq w(\gamma_i)\}$$
	where
	$$w(\gamma_i)=\mathop{\rm arcsinh} \left(\frac{1}{\sinh \frac{1}{2}\ell(\gamma_i)}\right)$$
	is the half width of the collar.
	
	Then the collars are pairwise disjoint for $i=1,...,m$. Each $\Cc(\gamma_i)$ is isomorphic to a cylinder $(\rho,t)\in [-w(\gamma_i),w(\gamma_i)] \times \mathbb S ^1$, where $\mathbb S ^1 = \R / \Z$, with the metric
\bear\label{collar-metric}	
ds^2=d\rho^2 + \ell(\gamma_i)^2 \cosh^2\rho dt^2.
\eear
	And for a point $(\rho,t)$, the point $(0,t)$ is its projection on the geodesic $\gamma_i$, the absolute value $\abs{\rho}$ is the distance to $\gamma_i$, $t$ is the coordinate on $\gamma_i \cong \mathbb S ^1$.
\end{theorem}
\noindent As the length $\ell(\gamma)$ of the central closed geodesic goes to $0$, the width \bear \label{wid-large} e^{w(\gamma)} \sim \frac{4}{\ell(\gamma)}\eear which tends to infinity. In this paper, we mainly deal with the case that $\ell(\gamma)$ is small and so $w(\gamma)$ is large.

Now we recall another version of the Collar Lemma which provides useful information on the injectivity radius at a point.
\begin{theorem}\cite[Theorem 4.1.6]{Buser10} \label{collar-inj}
Let $\beta_1, \cdots,\beta_k$ be the set of all simple closed geodesics of length $\leq 2\arcsinh(1)$ on a closed hyperbolic Riemann surface $X$ of genus $g$. Then $k\leq 3g-3$, and the followings hold.
\ben
\item The geodesics $\beta_1, \cdots, \beta_k$ are pairwise disjoint
\item $\inj_X(p)>\arcsinh(1)$ for any $p \in X \setminus (\cup_{i=1}^k \Cc(\beta_i))$.
\item If $p \in \Cc(\beta_i)$, and $d(p)=\dist(p,\p \Cc(\beta_i))$, then
\bear
\sinh(\inj_X(p))=\cosh\left(\frac{\ell(\beta_i)}{2}\right) \cosh\left(d(p))-\sinh(d(p)\right),\label{inj-1}\\
\quad \quad \quad   \sinh(\inj_X(p))=\sinh\left(\frac{\ell(\beta_i)}{2}\right) \cosh(\dist(p,\beta_i)) \ \emph{(see \cite[4.1.7]{Buser10})} \label{inj-2}
\eear
\een
\end{theorem}
\noindent By comparing the total area of all these standard collars $\Cc(\beta_i)'s$ and the total area $4\pi(g-1)$ of $X$, the set $X \setminus (\cup_{i=1}^k \Cc(\beta_i))$ is always non-empty. And for any point $q \in \p \Cc(\beta_i)$, by continuity or \eqref{inj-1} we know that $\inj_X(q) \geq \arcsinh(1)$. Now we study the injectivity radius along shortest geodesic loops (may not smooth at base points). First we consider the case that the base point is contained in a collar with a central closed geodesic of length $\leq 2\arcsinh(1)$.
\bpro \label{inj-short}
Let $X$ be a closed hyperbolic surface. For any $p \in X$ with $\inj_X(p)\leq \arcsinh(1)$, we let $\sigma:[0,2\inj_X(p)] \to X$ be a shortest nontrivial geodesic loop based at $p$. Then for any $s\in [0,2\inj_X(p)]$, we have
\[\left(\sqrt{2}-1\right)\inj_X(p)\leq \inj_X(\sigma(s))\leq \inj_X(p).\]
\epro  

\bp
Since $\inj_X(p)\leq \arcsinh(1)$, by Theorem \ref{collar} and \ref{collar-inj} one may assume that $\beta$ is the unique simple closed geodesic of length $\leq 2\arcsinh(1)$ such that $p \in \Cc(\beta)$. So the shortest geodesic loop $$\sigma([0,2\inj_X(p)])\subset \Cc(\beta)$$ otherwise it follows by Theorem \ref{collar-inj} that $\sigma([0,2\inj_X(p)])$ contains a point $q \notin \Cc(\beta)$ with $\inj_X(q)>\arcsinh(1)$ implying that $\ell(\sigma)>2\arcsinh(1)$ which is a contradiction. 

Now we first show the right hand side inequality: up to a conjugation one may lift $\beta$ onto the imaginary axis $\textbf{i}\cdot \R^+$ in the upper half plane $\mathbb{H}$, and the deck transformation corresponding to $\beta$ is $A(z)=e^{\ell(\beta)}z$. By Theorem \ref{collar} one may let $\tilde{p}$ be a lift of $p$ with $$\dist_{\bH}(\tilde{p},\textbf{i}\cdot \R^+)=\dist(p,\beta)\leq \arcsinh\left(\frac{1}{\sinh \frac{1}{2}\ell(\beta)}\right).$$ Then the lift $\tilde{\sigma}:[0,2\inj_X(p)] \to \mathbb{H}$ of $\sigma$ based at $\tilde{p}$ is the geodesic joining $\tilde{p}$ and $A(\tilde{p})$ in $\mathbb{H}$. By the convexity of distance functions on $\mathbb{H}$ we know that for any $s\in[0, 2\inj_X(p)]$,
\[\dist_{\bH}(\tilde{\sigma}(s), \textbf{i}\cdot \R^+)\leq \dist_{\bH}(\tilde{p},\textbf{i}\cdot \R^+)=\dist_{\bH}(A(\tilde{p}),\textbf{i}\cdot \R^+)\] 
which implies that
\bear\label{pi-decr}
\dist(\sigma(s),\beta)\leq \dist(p,\beta).
\eear
Then it follows by \eqref{inj-2} that $$\inj_X(\sigma(s))\leq \inj_X(p)$$ for all $s\in[0, 2\inj_X(p)]$.\\

Next we show the other side inequality. Let $s\in [0,2\inj_X(p)]$. Then
\bear\label{up-arcs1}
\dist(p,\sigma(s))\leq \inj_X(p)\leq \arcsinh(1).
\eear

We finish the proof by considering the following two cases.

\emph{Case $(1)$. $\dist(p,\beta)>\dist(p,\sigma(s))$.} Let $\pi(\sigma(s))\in \beta$ with $$\dist(\sigma(s),\pi(\sigma(s)))=\dist(\sigma(s),\beta).$$

\noindent By the triangle inequality we have 
\beqar
\dist(\sigma(s),\beta)&\geq& \dist(\pi(\sigma(s)),p)-\dist(p,\sigma(s))\\
&\geq& \dist(p,\beta)-\dist(p,\sigma(s))>0.
\eeqar

\noindent Then it follows by \eqref{inj-2} and \eqref{up-arcs1} that 
\beqar
\sinh(\inj_X(\sigma(s)))&=&\sinh\left(\frac{\ell(\beta)}{2}\right) \cosh(\dist(\sigma(s),\beta)) \nonumber \\
&\geq& \sinh\left(\frac{\ell(\beta)}{2}\right) \cosh\left(\dist(p,\beta)-\dist(p,\sigma(s)) \right) \\
&\geq & \left(\sinh\left(\frac{\ell(\beta)}{2}\right) \cosh(\dist(p,\beta))\right) e^{-\dist(p,\sigma(s))} \\
&\geq & \sinh (\inj_X(p))\cdot e^{-\arcsinh(1)} \\
&>&  \sinh \left(e^{-\arcsinh(1)}\cdot \inj_X(p)\right). 
\eeqar
Thus, we have
\bear \label{sqrt2-1}
\inj_X(\sigma(s))\geq e^{-\arcsinh(1)}\inj_X(p)=\left(\sqrt{2}-1\right)\inj_X(p).
\eear Which completes the proof for this case.

\emph{Case $(2)$. $\dist(p,\beta)\leq \dist(p,\sigma(s))$.} Then by \eqref{up-arcs1} we have 
\bear\label{pi-up}
\dist(p,\beta)\leq \dist(p,\sigma(s))\leq \arcsinh(1).
\eear
\noindent By \eqref{collar-metric} of Theorem \ref{collar}, we have that the closed curve (not a geodesic loop) based at $p$ with equidistance $\dist(p,\beta)$ to $\beta$ has length $\ell(\beta)\cosh(\dist(p,\beta))$. Which together with \eqref{pi-up} implies that
\[\inj_X(p)< \frac{\ell(\beta)}{2}\cosh(\arcsinh(1))=\frac{\sqrt{2}}{2}\ell(\beta).\]
\noindent Thus, we have
\bear \label{sqrt2-2}
\inj_X(\sigma(s))\geq \frac{\ell(\beta)}{2}> \frac{\sqrt{2}}{2}\inj_X(p).
\eear

Then the conclusion follows by \eqref{sqrt2-1} and \eqref{sqrt2-2} because $\frac{\sqrt{2}}{2}>(\sqrt{2}-1)$. 
\ep

Now we consider the case that the base point has injectivity radius larger than $\arcsinh(1)$. Let $\beta$ be a simple closed geodesic in $X$ of length $$\ell(\beta)\leq 2\arcsinh(1).$$ The boundary $\p \Cc(\beta)$ of the collar $\Cc(\beta)$ are two disjoint closed curves homotopic to $\beta$. By \eqref{collar-metric} we know that for each component $\beta'$ of $\p \Cc(\beta)$, the length $\ell(\beta')$ of $\beta'$ is
\beqar
\quad \ell(\beta')=\ell(\beta)\cosh(\omega(\beta))=\ell(\beta)\cdot \cosh\left( \arcsinh\left(\frac{1}{\sinh \frac{1}{2}\ell(\beta)}\right) \right).
\eeqar
A simple computation shows that
\bear \label{length-b}
\quad \ell(\beta')&=&\frac{\ell(\beta)}{\sinh \frac{1}{2}\ell(\beta)}\cdot \sqrt{1+\left(\sinh \frac{1}{2}\ell(\beta)\right)^2}\leq 2\sqrt{2}.
\eear
Now we are ready to state the result for the other case.
\bpro\label{inj-larg}
Let $X$ be a closed hyperbolic surface. For any $p \in X$ with $\inj_X(p)> \arcsinh(1)$, we let $\sigma:[0,2\inj_X(p)] \to X$ be a shortest nontrivial geodesic loop based at $p$. Then
\[\min_{s\in [0,2\inj_X(p)]} \inj_X(\sigma(s)) \geq \ln\left( e^{-\sqrt{2}}+\sqrt{e^{-2\sqrt{2}}+1} \right)\sim 0.2407.\]
\epro
\bp
If the geodesic loop $\sigma$ does not intersect with any standard collar with central closed geodesic of length less than $2\arcsinh(1)$, it follows by Theorem \ref{collar-inj} that 
\[\min_{s\in [0,2\inj_X(p)]} \inj_X(\sigma(s)) \geq \arcsinh(1)>\ln\left( e^{-\sqrt{2}}+\sqrt{e^{-2\sqrt{2}}+1} \right)\sim 0.2407.\]

Now we assume that $$\sigma \cap \Cc(\beta)\neq \emptyset$$ for some simple closed geodesic $\beta$ of length less than or equal to $2\arcsinh(1)$.

{\emph{Claim: $\forall q\in (\sigma \cap \Cc(\beta))$, $\dist(q,\partial \Cc(\beta))\leq \sqrt{2}$.}} 

If the claim above is true, then it follows by \eqref{inj-1} in Theorem \ref{collar-inj} that
\beqar
\sinh(\inj_X(q))&=&\cosh\left(\frac{\ell(\beta)}{2}\right) \cosh\left(\dist(q,\p \Cc(\beta)))-\sinh(\dist(q,\p \Cc(\beta))\right)\\
&\geq& \cosh\left(\dist(q,\p \Cc(\beta)))-\sinh(\dist(q,\p \Cc(\beta))\right)\\
&=& e^{-\dist(q,\p \Cc(\beta))}\\
&\geq& e^{-\sqrt{2}}
\eeqar
which implies that
\[\inj_X(q)\geq \ln\left( e^{-\sqrt{2}}+\sqrt{e^{-2\sqrt{2}}+1} \right)\sim 0.2407.\]
Since $\beta$ is an arbitrary closed geodesic of length less than or equal to $2\arcsinh(1)$ such that $\sigma \cap \Cc(\beta)\neq \emptyset$, by Theorem \ref{collar-inj} we know that
\[\min_{s\in [0,2\inj_X(p)]} \inj_X(\sigma(s)) \geq \min\left\{\arcsinh(1),\ 0.2407 \right\}=0.2407.\]
Now we prove the claim.

\noindent \emph{Proof of Claim.} Suppose for contradiction that $$\dist(q,\partial \Cc(\beta))> \sqrt{2} $$ for some $q \in (\sigma \cap \Cc(\beta))$. We let $t_1>0$ be the first time when $\sigma$ meets $\Cc(\beta)$, and $t_2>0$ be the last time when $\sigma$ meets $\Cc(\beta)$. That is,
\[t_1:=\min \left\{s \in [0,2\inj_X(p)]; \ \sigma(s)\in \Cc(\beta)\right\}\] 
and
\[t_2:=\max \left\{s \in [0,2\inj_X(p)]; \ \sigma(t)\in \Cc(\beta)\right\}.\] 
Since $\inj_X(\sigma(0))=\inj_X(\sigma(2\inj_X(p)))=\inj_X(p)>\arcsinh(1)$, we have that $t_1>0$ and $t_2<2\inj_X(p)$. Clearly we have $\sigma(t_1)\cup\sigma(t_2)\in \p \Cc(\beta)$. If $\sigma(t_1)$ and $\sigma(t_2)$ are on the same component of the boundary $\p \Cc(\beta)$, then we have
\[t_2-t_1\geq \dist(\sigma(t_1),q)+\dist(q,\sigma(t_2))>2\sqrt{2}.\]   
If $\sigma(t_1)$ and $\sigma(t_2)$ are on the different components of the boundary $\p \Cc(\beta)$, by symmetry of the standard collar $\Cc(\beta)$ we also have
\[t_2-t_1\geq 2\dist (q,\p \mathcal{C}(\beta))> 2\sqrt{2}.\]
So we always have 
\bear \label{gap-2}
t_2-t_1>2\sqrt{2}.
\eear

Now we finish the argument by considering the following two cases.

\emph{Case $(1)$. $t_1=\ell(\sigma([0,t_1]))\leq \ell(\sigma([t_2,2\inj_X(p)]))=2\inj_X(p)-t_2$.} Let $\beta'$ be the component of the boundary $\p \Cc(\beta)$ with $\sigma(t_1)\in \beta'$, and we parametrize $\beta'$ such that $\beta'(0)=\beta'(\ell(\beta'))=\sigma(t_1)$. Consider the closed curve $\sigma'$ based at $p$ as following.
\beqar
\sigma'(s):=
\begin{cases}
\sigma(s), \quad \quad \quad \ \emph{$s\in [0,t_1]$}; \\

\beta'(s-t_1),  \quad \emph{$s\in [t_1, t_1+\ell(\beta')]$};\\

\sigma\left( 2t_1+\ell(\beta')-s\right), \quad \emph{$s\in [t_1+\ell(\beta'), 2t_1+\ell(\beta')]$}.
\end{cases}
\eeqar
The closed curve $\sigma'$ is freely homotopic to $\beta$. So $\sigma'$ is nontrivial. By \eqref{length-b} and \eqref{gap-2}, the length of $\sigma'$ satisfies 
\beqar
\ell(\sigma')&=&2t_1+\ell(\beta') \\
&\leq& t_1+(2\inj_X(p)-t_2)+2\sqrt{2} \\
&<& 2\inj_X(p)
\eeqar
which is a contradiction since $\sigma$ is a shortest nontrivial geodesic loop based at $p$.

\emph{Case $(2)$. $t_1=\ell(\sigma([0,t_1]))\geq \ell(\sigma([t_2,2\inj_X(p)]))=2\inj_X(p)-t_2$.} Let $\beta''$ be the component of the boundary $\p \Cc(\beta)$ with $\sigma(t_2)\in \beta''$. Similarly, we consider the closed curve $\sigma''$ based at $p$ which is defined as $$\sigma''=\sigma([t_2,2\inj_X(p)])\cup \beta'' \cup \sigma([t_2,2\inj_X(p)]).$$ Then the closed curve $\sigma''$ is freely homotopic to $\beta$. So $\sigma''$ is nontrivial. By \eqref{length-b} and \eqref{gap-2}, the length of $\sigma''$ satisfies 
\beqar
\ell(\sigma'')&=&2(2\inj_X(p)-t_2)+\ell(\beta'') \\
&\leq& t_1+(2\inj_X(p)-t_2)+2\sqrt{2} \\
&<& 2\inj_X(p)
\eeqar
which is a contradiction since $\sigma''$ is a nontrivial closed loop based at $p$.

The proof is complete.
\ep

\section{One useful inequality}\label{sec-ine}
In this section we prove the following property which is crucial in the proof of Theorem \ref{mt}.
\bpro\label{neck-ineq}
Let $X$ be a hyperbolic surface. For any $p \in X$ we let $\sigma:[0,2\inj_X(p)] \to X$ be a shortest nontrivial geodesic loop based at $p$. Assume that $$\inf_{s\in [0,2\inj_X(p)]}\inj_X(\sigma(s))\geq 2\eps_0$$ for some uniform constant $\eps_0>0$. Then for any function $f\geq 0$ on $X$, we have
\[\int_0^{2\inj_X(p)} \left(\int_{B(\sigma(s);\eps_0)}f\dArea\right) ds\leq 12 \eps_0 \int_{\N_{\eps_0}(\sigma)} f \dArea.\]
\noindent Where $B(\sigma(s);\eps_0)=\{q\in X; \ \dist(q, \sigma(s))<\eps_0\}$ and $\N_{\eps_0}(\sigma)$ is the $\eps_0$-\nbhd \ of $\sigma$, i.e., 
\[\N_{\eps_0}(\sigma)=\{z\in X; \ \dist\left(x, \sigma([0,2\inj_X(p)])\right)< \eps_0\}\]
\epro

We split the proof into several parts.

First since $\sigma:[0,2\inj_X(p)] \to X$ is a shortest nontrivial geodesic loop based at $p$, it is known that both the two restrictions $\sigma:[0,\inj_X(p)] \to X$ and $\sigma:[\inj_X(p),2\inj_X(p)]\to X$ are minimizing geodesics. For any $s\in (0,2\inj_X(p))$, we let $\vec{n}(s)$ be an unit normal vector of $\sigma$ at $\sigma(s)$. Consider the foliation $\left\{\exp_{\sigma(s)}(t\cdot \vec{n}(s))\right\}_{t\in (-\eps_0,\eps_0), \ t\in (0,2\inj_X(p))}$ along $\sigma$ where $\exp_{\sigma(s)}(\cdot)$ is the standard exponential map at $\sigma(s)$. Set
\[m=\left[\frac{\inj_X(p)}{\eps_0}\right]\]
to be the largest integer of the number $\frac{\inj_X(p)}{\eps_0}$. 

Now we assume that $$m\geq 2.$$ 

\begin{figure}[htbp]
   \centering
   \includegraphics[width=3.5in]{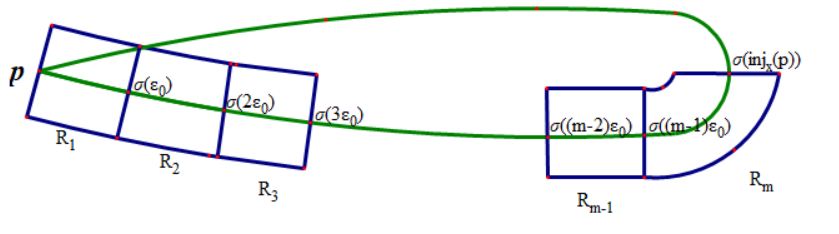} 
   \caption{$R_i$'s and $R_m$}
   \label{loop}
\end{figure}

Set
\ben
\item for $1\leq i \leq m-1$, $$R_i=\bigcup_{s\in \left[(i-1)\eps_0, i\eps_0\right)} \ \bigcup_{t\in \left(-\eps_0,\eps_0\right)}\exp_{\sigma(s)}(t\cdot \vec{n}(s))$$ 
 and $$R_m=\bigcup_{s\in \left[(m-1)\eps_0, \inj_X(p)\right]}\ \bigcup_{t\in \left(-\eps_0,\eps_0\right)}\exp_{\sigma(s)}(t\cdot \vec{n}(s))$$
(see Figure \ref{loop}).

\item for $1\leq i \leq m-1$, $$R_{i}^{'}=\bigcup_{s\in \left(2\inj_X(p)-i\cdot\eps_0,2\inj_X(p)-(i-1)\eps_0\right]}\ \bigcup_{t\in \left(-\eps_0,\eps_0\right)}\exp_{\sigma(s)}(t\cdot \vec{n}(s))$$ 
and $$R_{m}^{'}=\bigcup_{s\in \left[\inj_X(p),2\inj_X(p)-(m-1)\eps_0\right]} \ \bigcup_{t\in \left(-\eps_0,\eps_0\right)}\exp_{\sigma(s)}(t\cdot \vec{n}(s)).$$
\een

\bl \label{disjoint} With the notations above, 
\ben
\item  $\bigcup_{i=1}^m \left(R_i \cup R_{i}^{'}\right)\subseteq \N_{\eps_0}(\sigma);$

\item for all $1\leq i \neq j \leq m$, we have
\[R_i \cap R_j=\emptyset \quad \emph{and} \quad R_{i}^{'} \cap R_{j}^{'}=\emptyset.\]
\een
\el
\bp
Part (1) is clear.

For Part $(2)$, we first prove $$R_i \cap R_j =\emptyset \ \emph{for} \ 1\leq i \neq j \leq m.$$ Suppose for contradiction that there would exist a point $q \in R_i \cap R_j$ for some $i\neq j \in [1,m]$. Recall that $\sigma:[0,\inj_X(p)]\to X$ is a minimizing geodesic. By construction one may assume that $q_i\in R_i$ and $q_j\in R_j$ such that
\ben
\item[(i)] the geodesic triangle $\Delta(p_i,p_j,q)$ with vertices $p_i,p_j$ and $q$ has at least two interior $\frac{\pi}{2}$-angles;
\item[(ii)] $\max\{\dist(p_i,q),\dist(p_j,q)\}< \eps_0$.
\een
It follows by the triangle inequality and Part (ii) above that $$\dist(p_i,p_j)< 2\eps_0$$ which implies that the geodesic triangle $\Delta(p_i,p_j,q)\subset B(p_i;2\eps_0)$. Recall that $\inj_X(p_i)\geq 2\eps_0$. So $\Delta(p_i,p_j,q)\subset B(p_i;2\eps_0)$ bounds a disk. By the Gauss-Bonnet formula \cite{Buser10} we know that the total interior angle of $\Delta(p_i,p_j,q)$ is less than $\pi$, which contradicts Part $(i)$. 

The proof for $ R_{i}^{'} \cap R_{j}^{'}=\emptyset$ is the same as above.
\ep

\br
It is possible that for all $1\leq i \leq m$, $R_i \cap R_{i}^{'}\neq \emptyset$. 
\er

\bl \label{ends} With the notations above, 
\ben
\item $\int_0^{\eps_0} \left(\int_{B(\sigma(s);\eps_0)}f\dArea\right) ds\leq \eps_0 \int_{\N_{\eps_0}(\sigma)} f \dArea.$

\item $\int_{(m-1)\eps_0}^{\inj_X(p)} \left(\int_{B(\sigma(s);\eps_0)}f\dArea\right) ds\leq 2 \eps_0 \int_{\N_{\eps_0}(\sigma)} f \dArea.$
\een
\el
\bp
It follows by the triangle inequality that
\[\cup_{0\leq s \leq \eps_0}B(\sigma(s);\eps_0) \subset \N_{\eps_0}(\sigma) \ \emph{and} \ \cup_{(m-1)\eps_0\leq s \leq \inj_X(p)}B(\sigma(s);\eps_0) \subset \N_{\eps_0}(\sigma)\]
which implies $(1)$, and 
\beqar
\int_{(m-1)\eps_0}^{\inj_X(p)} \left(\int_{B(\sigma(s);\eps_0)}f\dArea\right) &\leq& (\inj_X(p)-(m-1)\eps_0)  \int_{\N_{\eps_0}(\sigma)} f \dArea \\
&\leq& 2 \eps_0 \int_{\N_{\eps_0}(\sigma)} f\dArea.
\eeqar
The proof is complete.
\ep

Now we are ready to prove Proposition \ref{neck-ineq}.
\bp [Proof of Proposition \ref{neck-ineq}]
If $$m=\left[\frac{\inj_X(p)}{\eps_0}\right]\leq 5,$$ then $2\inj_X(p)\leq 12 \eps_0$. Since $\cup_{0\leq s \leq 2\inj_X(p)}B(\sigma(s);\eps_0) \subset \N_{\eps_0}(\sigma),$
we have
\[\int_0^{2\inj_X(p)} \left(\int_{B(\sigma(s);\eps_0)}f\dArea\right) ds\leq 12 \eps_0 \int_{\N_{\eps_0}(\sigma)} f \dArea\]
which completes the proof.

Now we always assume that $$m=\left[\frac{\inj_X(p)}{\eps_0}\right]\geq 6.$$
From Lemma \ref{ends} we have
\bear \label{n-e-1}
&& \quad \int_0^{\inj_X(p)} \left(\int_{B(\sigma(s);\eps_0)}f\dArea\right) ds=\int_0^{\eps_0} \left(\int_{B(\sigma(s);\eps_0)}f\dArea\right)ds+ \\
&& \sum_{i=1}^{m-2}\int_{i \eps_0}^{(i+1)\eps_0} \left(\int_{B(\sigma(s);\eps_0)}f\dArea\right) ds+\int_{(m-1)\eps_0}^{\inj_X(p)} \left(\int_{B(\sigma(s);\eps_0)}f\dArea\right) ds \nonumber \\
&&  \leq 3\eps_0\int_{\N_{\eps_0}(\sigma)} f \dArea+ \sum_{i=1}^{m-2}\int_{i \eps_0}^{(i+1)\eps_0} \left(\int_{B(\sigma(s);\eps_0)}f\dArea\right)ds. \nonumber
\eear

\noindent For each $1\leq i \leq m-2$, by the triangle inequality we know that
\[ \bigcup_{s=i\eps_0}^{(i+1)\eps_0}B(\sigma(s);\eps_0)\subset (R_i\cup R_{i+1}\cup R_{i+2}).\]
Thus, by Lemma \ref{disjoint} we have
\beqar
&& \sum_{i=1}^{m-2}\int_{i \eps_0}^{(i+1)\eps_0} \left(\int_{B(\sigma(s);\eps_0)}f\dArea\right)ds \leq  \sum_{i=1}^{m-2} \eps_0 \left(\int_{R_i\cup R_{i+1}\cup R_{i+2}}f\dArea\right) \\
&&=\eps_0  \sum_{i=1}^{m-2} \left(\int_{R_i}f\dArea+\int_{R_{i+1}}f\dArea+\int_{R_{i+2}}f\dArea\right) \nonumber \\
&&\leq 3 \eps_0 \int_{\bigcup_{i=1}^m R_i}f\dArea \nonumber \\
&& \leq 3\eps_0 \int_{\N(\sigma)}f\dArea \nonumber
\eeqar
which together with \eqref{n-e-1} implies that
\bear \label{n-e-2}
\quad \int_0^{\inj_X(p)} \left(\int_{B(\sigma(s);\eps_0)}f\dArea\right) ds\leq 6\eps_0\int_{\N(\sigma)}f\dArea.
\eear

Restricted on the geodesic $\sigma([\inj_X(p),2\inj_X(p)])$ and replace each $R_i$ by $R_i'$ , one may also apply the same argument above to get
\bear \label{n-e-3}
\quad \int_{\inj_X(p)}^{2\inj_X(p)} \left(\int_{B(\sigma(s);\eps_0)}f\dArea\right) ds\leq 6\eps_0\int_{\N(\sigma)}f\dArea.
\eear

Then the conclusion follows by \eqref{n-e-2} and \eqref{n-e-3}.
\ep

\br
The two main ingredients in the proof above are:
\ben
\item the two restrictions $\sigma:[0,\inj_X(p)] \to X$ and $\sigma:[\inj_X(p),2\inj_X(p)]\to X$ are minimizing geodesics;

\item the total interior angle of any geodesic triangle, which bounds a disk, is less than or equal to $\pi$. 
\een
Actually in the argument above, if one replaces the foliation $$\left\{\exp_{\sigma(s)}(t\cdot \vec{n}(s))\right\}_{t\in (-\eps_0,\eps_0), t\in (0,2\inj_X(p))}$$ along $\sigma$ by 
\[\left \{\exp_{\sigma(s)}(t\cdot \vec{w}(s))\right \}_{t\in (-\eps_0,\eps_0), \ t\in (0,2\inj_X(p)), \ \emph{$\vec{w}(s)$ is unit and normal to $\sigma'(s)$ }},\]
then it follows by the same argument above that one may generalize Proposition \ref{neck-ineq} to higher dimensions. More precisely, we have
\bpro \label{neck-ineq-h}
Let $M$ be a complete Riemannian manifold of nonpositive curvature. For any $p \in M$ we let $\sigma:[0,2\inj_M(p)] \to M$ be a shortest nontrivial geodesic loop based at $p$. Assume that $$\inf_{s\in [0,2\inj_M(p)]}\inj_X(\sigma(s))\geq 2\eps_0$$ for some uniform constant $\eps_0>0$. Then for any function $f\geq 0$ on $M$, we have
\[\int_0^{2\inj_M(p)} \left(\int_{B(\sigma(s);\eps_0)}f\dVol\right) ds\leq 12 \eps_0 \int_{\N_{\eps_0}(\sigma)} f \dVol.\]
\noindent Where $\N_{\eps_0}(\sigma)$ is the $\eps_0$-\nbhd \ of $\sigma$, i.e.,
\[\N_{\eps_0}(\sigma)=\{z\in X; \ \dist\left(x, \sigma([0,2\inj_M(p)])\right)< \eps_0\}.\]
\epro
\er

\section{Proof of Theorem \ref{mt}}\label{sec-mt}
In this section we prove Theorem \ref{mt}.

Let $c:[0,T_0]\to \sM_{-1}$ be a smooth horizontal path of \WP \ arc-length parameter where $T_0>0$ is a constant. For any $t\in (0,T_0)$, here one may view the tangent vector $c'(t)$ as a harmonic Beltrami differential on the hyperbolic surface $c(t) \in \sM_{-1}$. Since $t$ is an arc-length parameter, $$\int_{c(t)}|c'(t)(z)|^2\dArea(z) \equiv 1$$ for all $t\in(0,T_0)$. It is known from \cite{RT18} that this function $\inj_{c(t)}(p)$ is differentiable almost everywhere in $(0,T_0)$ and the Fundamental Theorem of Calculus  holds for $\inj_{c(t)}(p)$ along $c$. Now we recall the following two lemmas from \cite{RT18}, whose proofs are outlined here for completeness.

\bl[Rupflin-Topping,\ {\cite[Lemma 2.1]{RT18}}]\label{inj-lip}
Let $c:[0,T_0]\to \sM_{-1}$ be a a smooth horizontal path of \WP \ arc-length parameter where $T_0>0$ is a fixed constant. Then for any $p \in S_g$, the function $\inj_{c(t)}(p):[0,T_0]\to \R^{>0}$ is locally Lipschitz.
\el

\bp
For completeness we outline the proof here. One may see the proof of \cite[Lemma 2.1]{RT18} for more details. Let $I\subset [0,T_0]$ be any sub-interval. Since $c(t)$ is smooth, for any $t_1,t_2\in I$ there exists a constant $C>0$ such that for all $t_1,t_2,t\in I$,
\[||c(t_1)-c(t_2)||_{c(t)}\leq C\cdot |t_1-t_2|\]
where $||c(t_1)-c(t_2)||_{c(t)}$ is the norm of the difference of two hyperbolic metrics $c(t_1)$ and $c(t_2)$ at $c(t)$. Let 
$\sigma_i (i=1,2):[0,2\inj_{c(t_i)}(p)]\to c(t_i)$ be a shortest geodesic loop based at $p$ respectively. Without loss of generality one may assume that $\inj_{c(t_1)}(p)\geq \inj_{c(t_2)}(p)$. Since,
\beqar
2\inj_{c(t_1)}(p)&\leq& \int_{0}^{2\inj_{c(t_2)}(p)}\sqrt{\left<\sigma_2'(s), \sigma_2'(s)\right>_{c(t_1)}}ds\\
&\leq & \frac{1}{2} \int_{0}^{2\inj_{c(t_2)}(p)}\left(\left<\sigma_2'(s), \sigma_2'(s)\right>_{c(t_1)}+1  \right)ds.
\eeqar
Thus, we have
\beqar
&&2(\inj_{c(t_1)}(p)-\inj_{c(t_2)}(p))\leq\frac{1}{2} \int_{0}^{2\inj_{c(t_2)}(p)}\left(\left<\sigma_2'(s), \sigma_2'(s)\right>_{c(t_1)}-1  \right)ds \\
&&= \frac{1}{2} \int_{0}^{2\inj_{c(t_2)}(p)}\left(\left<\sigma_2'(s), \sigma_2'(s)\right>_{c(t_1)}-\left<\sigma_2'(s), \sigma_2'(s)\right>_{c(t_2)}  \right)ds\\
&&\leq \frac{1}{2} \int_{0}^{2\inj_{c(t_2)}(p)}||c(t_2)-c(t_1)||_{c(t_2)}ds\\
&& \leq C \inj_{c(t_2)}(p)|t_1-t_2|\\
&& \leq C \ln(4g-2)|t_1-t_2|
\eeqar
which completes the proof.
\ep

\bl[Rupflin-Topping, \ {\cite[Lemma 2.5]{RT18}}]\label{inj-de}
Let $c:[0,T_0]\to \sM_{-1}$ be a smooth horizontal path of \WP \ arc-length parameter where $T_0>0$ is a fixed constant. For any $p \in S_g$, and suppose that the function $\inj_{c(t)}(p):[0,T_0]\to \R^{>0}$ is differentiable at $t=t_0 \in (0,T_0)$. Then for any shortest geodesic loop $\sigma:[0,2\inj_{c(t_0)}(p)]\to c(t_0)$ based at $p$,
\[\frac{d}{dt}\inj_{c(t)}(p)\big|_{t=t_0}= \frac{1}{4}\int_{0}^{2\inj_{c(t_0)}(p)} c'(t_0)(\sigma'(s),\sigma'(s))ds.\]
\el
\bp
For completeness we also provide the proof here. Set \[f(t):=\frac{1}{2}\int_{0}^{2\inj_{c(t_0)}(p)}\sqrt{\left<\sigma'(s), \sigma'(s)\right>_{c(t)}}ds-\inj_{c(t)}(p).\]
So we have $f(t)\geq 0$ and $f(t_0)=0$. Then the conclusion follows by that $f'(t_0)=0$.
\ep

Now we prove the key estimation in the proof of Theorem \ref{mt}.
\bpro \label{key}
Let $c:[0,T_0]\to \sM_{-1}$ be a smooth horizontal path of \WP \ arc-length parameter where $T_0>0$ is a fixed constant. For any $p \in S_g$, and suppose that the function $\inj_{c(t)}(p):[0,T_0]\to \R^{>0}$ is differentiable at $t=t_0 \in (0,T_0)$. Then we have
\[\left|\frac{d}{dt}\sqrt{\inj_{c(t)}(p)}\bigg|_{t=t_0}\right|\leq 0.3884.\]
\epro
\begin{rem*}
It was shown in \cite[Lemma 2.2]{RT18} by Rupflin and Topping that $\left|\frac{d}{dt}\sqrt{\inj_{c(t)}(p)}\bigg|_{t=t_0}\right|\leq K_g$ where $K_g>0$ is a constant depending on $g$. Our essential improvement for Proposition \ref{key} is that this constant $K_g$ can be chosen to be uniform.
\end{rem*}

\bp [Proof of Proposition \ref{key}]
We follow the same idea as the proof of \cite[Lemma 2.2]{RT18}, and prove it by two cases. Let $\sigma:[0,2\inj_{c(t_0)}(p)]\to c(t_0)$ be a shortest geodesic loop based at $p$. Thus, it follows by Lemma \ref{inj-de} that
\bear\label{k-eq-2}
\left|\frac{d}{dt}\inj_{c(t)}(p)\bigg|_{t=t_0}\right|&\leq& \frac{1}{4}\int_{0}^{2\inj_{c(t_0)}(p)} \left|c'(t_0)(\sigma'(s),\sigma'(s))\right|ds \\
&\leq & \frac{1}{4}\int_{0}^{2\inj_{c(t_0)}(p)} \left|c'(t_0)(\sigma(s))\right|ds. \nonumber
\eear
Either $(1). \inj_{c(t_0)}(p)\leq \arcsinh(1)$ or $(2).\inj_{c(t_0)}(p)>\arcsinh(1)$. We finish the proof by considering these two cases.

\emph{Case $(1)$. $\inj_{c(t_0)}(p)\leq \arcsinh(1)$.} By Proposition \ref{inj-short} we know that for any $s\in [0,2\inj_{c(t_0)}(p)]$, 
\bear\label{k-eq-1}
\quad \left(\sqrt{2}-1\right)\inj_{c(t_0)}(p)\leq \inj_{c(t_0)}(\sigma(s))\leq \inj_{c(t_0)}(p)\leq \arcsinh(1).
\eear

\noindent Then one may apply Proposition \ref{bound-BW} to get that for any $s\in [0,2\inj_{c(t_0)}(p)]$, 
\bear\label{k-eq-3}
 \left|c'(t_0)(\sigma(s))\right| &\leq& \sqrt{\frac{\int_{c(t_0)}|c'(t_0)(z)|^2\dArea(z)}{\inj_{c(t_0)}(\sigma(s))}}\\
&=&\frac{1}{\sqrt{\inj_{c(t_0)}(\sigma(s))}} \nonumber
\eear
where in the last inequality we apply the fact that $t$ is an arc-length parameter. Thus, by \eqref{k-eq-2}, \eqref{k-eq-1} and \eqref{k-eq-3} we have
\beqar
\left|\frac{d}{dt}\inj_{c(t)}(p)\bigg|_{t=t_0}\right| &\leq & \frac{1}{4}\int_{0}^{2\inj_{c(t_0)}(p)}  \frac{1}{\sqrt{\inj_{c(t_0)}(\sigma(s))}}ds \\
&\leq & \frac{\sqrt{\inj_{c(t_0)}(p)}}{2\sqrt{\sqrt{2}-1}}
\eeqar
which implies that
\bear \label{g-bound-short}
\left|\frac{d}{dt}\sqrt{\inj_{c(t)}(p)}\bigg|_{t=t_0}\right|\leq \frac{1}{4\sqrt{\sqrt{2}-1}}\sim 0.3884.
\eear

\emph{Case $(2)$. $\inj_{c(t_0)}(p)>\arcsinh(1)$.} First by Proposition \ref{inj-larg} we know that
\bear \label{k-eq-low}
\min_{s\in [0,2\inj_X(p)]} \inj_X(\sigma(s)) \geq  0.2407.
\eear
Let $$0<r_0\leq \frac{0.2407}{2}\sim 0.1203.$$ One may apply Proposition \ref{bound-teo} to get that for any $s\in [0,2\inj_{c(t_0)}(p)]$, 
\bear\label{k-eq-4}
 \quad \left|c'(t_0)(\sigma(s))\right| \leq C(r_0) \sqrt{\int_{B_{c(t_0)}(\sigma(s);r_0)}|c'(t_0)(z)|^2 \dArea(z)}
\eear
where $B_{c(t_0)}(\sigma(s);r_0)=\left\{q\in c(t_0); \ \dist(q, \sigma(s))<r_0 \right\}$. Then by \eqref{k-eq-2} and \eqref{k-eq-4} we have
\small{
\bear
&&\left|\frac{d}{dt}\inj_{c(t)}(p)\bigg|_{t=t_0}\right| \leq \frac{C(r_0)}{4}\int_{0}^{2\inj_{c(t_0)}(p)} \sqrt{ \int_{B_{c(t_0)}(\sigma(s);r_0)}|c'(t_0)(z)|^2 \dArea(z)}ds \nonumber \\
&&\leq  \frac{C(r_0)}{4} \sqrt{2\inj_{c(t_0)}(p)}\sqrt{ \int_{0}^{2\inj_{c(t_0)}(p)}\left( \int_{B_{c(t_0)}(\sigma(s);r_0)}|c'(t_0)(z)|^2 \dArea(z)\right) ds} \nonumber
\eear
}
\noindent where in the last inequality we apply the Cauchy-Schwarz inequality. In light of \eqref{k-eq-low}, we now apply Proposition \ref{neck-ineq}. Let $X=c(t_0)$, $\eps_0=r_0$ and $f(z)=|c'(t_0)(z)|^2\geq 0$ in Proposition \ref{neck-ineq}, then it follows by Proposition \ref{neck-ineq} that
\bear
\left|\frac{d}{dt}\inj_{c(t)}(p)\bigg|_{t=t_0}\right| &\leq & \frac{\sqrt{6r_0}}{2}C(r_0)  \sqrt{\inj_{c(t_0)}(p)} \sqrt{\int_{\N_{r_0}(\sigma)} |c'(t_0)(z)|^2 \dArea(z)} \nonumber\\
&\leq&\frac{\sqrt{6r_0}}{2}C(r_0) \sqrt{\inj_{c(t_0)}(p)} \nonumber
\eear
where in the last inequality we apply that $c(t)$ is of \WP \ arc-length parameter. Thus, we have
\bear\label{g-bound-larg-0}
\left|\frac{d}{dt}\sqrt{\inj_{c(t)}(p)}\bigg|_{t=t_0}\right|\leq \frac{\sqrt{6}}{4}C(r_0)\sqrt{r_0}.
\eear
Recall that $C(r) = \frac{1}{\sqrt{\pi r}} + O(1)$ as $r\to 0$. Let $r_0 \to 0$ in \eqref{g-bound-larg-0} we get
\bear\label{g-bound-larg}
\left|\frac{d}{dt}\sqrt{\inj_{c(t)}(p)}\bigg|_{t=t_0}\right|\leq \frac{\sqrt{6}}{4\sqrt{\pi}}\sim0.3454.
\eear

Then the conclusion follows by \eqref{g-bound-short} and \eqref{g-bound-larg}.
\ep

Now we are ready to prove Theorem \ref{mt}.
\bt[= Theorem \ref{mt}]
Fix a point $p \in S_g$ $(g\geq2)$. Then for any $X, Y \in \sT_g$, 
$$\left|\sqrt{\inj_X(p)}-\sqrt{\inj_Y(p)}\right|\leq  0.3884 \dist_{wp}(X,Y)$$
where $\dist_{wp}$ is the \WP \ distance.
\et

\bp[Proof of Theorem \ref{mt}]
For any $X,Y \in \Teich(S_g)$, by Wolpert \cite{Wolpert87} there exists a unique \WP \ geodesic $\underline{c}:[0,\dist_{wp}(X,Y)]\to \Teich(S_g)$ of arc-length parameter such that $\underline{c}(0)=X$ and $\underline{c}(\dist_{wp}(X,Y))=Y$. We lift $\underline{c}$ onto a horizontal curve $c:[0,\dist_{wp}(X,Y)]\to \sM_{-1}$ such that $\pi(c(t))=\underline{c}(t)$ for all $t\in [0, \dist_{wp}(X,Y)$. By Lemma \ref{inj-lip} we know that the injectivity radius  $\inj_{c(t)}(p):[0,\dist_{wp}(X,Y)] \to \R^{>0}$ is locally Lipschitz. Then we apply the Fundamental Theorem of Calculus and Proposition \ref{key} to get
\beqar
|\sqrt{\inj_{c(0)}(p)}-\sqrt{\inj_{c(1)}(p)}|&=&\left|\int_{0}^{\dist_{wp}(X,Y)}\frac{d}{dt}\left(\sqrt{\inj_{c(t)}(p)}\right)dt\right| \\
&\leq &\int_{0}^{\dist_{wp}(X,Y)}\left|\frac{d}{dt}\left(\sqrt{\inj_{c(t)}(p)}\right)\right|dt \\
&\leq&  0.3884 \dist_{wp}(X,Y)
\eeqar
which implies the conclusion by letting $c$ runs over all such horizontal lifts of the \WP \ geodesic $\underline{c}$ joining $X$ and $Y$.
\ep

If we only restrict the proof of Theorem \ref{mt} on the Case $(2)$ in the proof of Proposition \ref{key}, one may get 
\begin{corollary} \label{mt-larg}
Fix a point $p \in S_g$ $(g\geq2)$ and let $X,Y\in \sT_g$. Assume that for some horizontal curve $c\subset \sM_{-1}$ such that the projection $\pi(c)$ is the \WP \ geodesic joining $X$ and $Y$, and 
\[\min_{t\in [0,\dist_{wp}(X,Y)]}\inj_{c(t)}(p)>\arcsinh(1),\]
then we have
\[\left|\sqrt{\inj_X(p)}-\sqrt{\inj_Y(p)}\right|\leq  0.3454 \dist_{wp}(X,Y).\]
\end{corollary}

\section{The geometry of $\sM_g$ for large genus}\label{mg-large}
In this section we study the asymptotic geometry of $\sM_g$ for large genus.

Before proving Theorem \ref{mt-3}, we recall several things in \cite{W-inradius}. 
Let $\alpha \subset X$ be a simple closed geodesic. Up to a conjugation one may lift $\alpha$ to the imaginary axis $\textbf{i}\mathbb{R}^{+}$ in the upper half plane $\bH$. A special case of Riera's formula \cite[Theorem 2]{Rie05} says that 
\begin{equation}\label{Rie-f}
\langle\,\grad \ell_{\alpha},\grad \ell_{\alpha}\rangle_{wp}(X)=\frac{2}{\pi}(\ell_{\alpha}+\displaystyle \sum_{B\in \left\{\left<A\right>\backslash\Gamma/\left<A\right>-id\right\}} (u \ln{\frac{u+1}{u-1}}-2))
\end{equation}
where $u=\cosh{(\dist_{\mathbb{H}}(\textbf{i}\mathbb{R}^{+}, B\circ \textbf{i}\mathbb{R}^{+}))}$ and the double-coset of the identity element is omitted from the sum.

From now on, we always assume that $\alpha$ is a systolic curve of $X$ with large length, more precisely,
\bear \label{ass-8}
\ell_{\alpha}(X)=\lsys\geq 8.
\eear

As in \cite[Page 1327]{W-inradius}, we know that for any $B\in \left\{\left<A\right>\backslash\Gamma/\left<A\right>-id\right\}$ there exists a unique point $p_B \in B\circ(\textbf{i}\mathbb{R}^{+})$ such that 
\[\dist_{\mathbb{H}}(p_B, \textbf{i}\mathbb{R}^{+})=\dist_{\mathbb{H}}(B\circ (\textbf{i}\mathbb{R}^{+}),\textbf{i}\mathbb{R}^{+}).\]
\noindent By \cite[Lemma 4.6]{W-inradius} and \cite[Lemma 4.8]{W-inradius}, one may choose a representative $B'\in \left<A\right>\backslash \Gamma-id$ for $B$ such that
\ben
\item $\dist_{\bH}(p_{B'},\textbf{i}\mathbb{R}^{+})\geq \frac{\lsys}{4}\geq 2$;

\item $1\leq r_{B'}\leq e^{\lsys}$
\een 
where $p_{B'}=(r_{B'},\theta_{B'})$ in polar coordinate be the nearest projection point on $ B'\circ (\textbf{i}\mathbb{R}^{+})$ from $\textbf{i}\mathbb{R}^{+}$. For $z=(r,\theta) \in \mathbb{H}$ given in polar coordinate where $\theta \in (0,\pi)$, the hyperbolic distance between $z$ and the imaginary axis $\textbf{i}\mathbb{R^+}$ is
\begin{eqnarray}\label{i-dis}
\dist_{\mathbb{H}}(z, \textbf{i}\mathbb{R^+})=\ln|\csc{\theta}+|\cot{\theta}||.
\end{eqnarray}
Which implies
\begin{eqnarray}\label{i-exp}
e^{-2\dist_{\mathbb{H}}(z, \textbf{i}\mathbb{R^+})}\leq \sin^2{\theta}=\frac{\Im^2(z)}{|z|^2}\leq 4e^{-2\dist_{\mathbb{H}}(z, \textbf{i}\mathbb{R^+})}.
\end{eqnarray}

Now we consider the geodesic balls $\{B(p_B; 1)\}_{B\in \left<A \right> \backslash \Gamma-id}$ of radius $1$ in $\bH$. 
\bl\label{unit-disjoint}
For any $B_1\neq B_2 \in  \left< A \right> \backslash \Gamma-id$,
\[B(p_{B_1};1)\cap B(p_{B_2};1)=\emptyset.\]
\el
\bp
It follows by the triangle inequality and \cite[Lemma 4.6]{W-inradius}.
\ep

\bl\label{range-union}
{\small{\[\bigcup_{B \in  \left<A \right> \backslash \Gamma-id}B(p_B; 1) \subset \left\{(r,\theta)\in \bH; \ e^{-1}\leq r \leq e^{\lsys+1} \ \emph{and} \ \sin(\theta)\leq 2e^{-\frac{\lsys}{8}} \right\}.\]}}
\el

\bp
For any $z=(r,\theta)\in (p_B; 1)$ where $B \in  \{\left<A \right> \backslash \Gamma-id\}$ is arbitrary, since $1\leq r_{B}\leq e^{\lsys}$, by the triangle inequality we clearly have \[ e^{-1}\leq r \leq e^{\lsys+1}.\]

Now we control the angle $\theta$. Since  $\dist_{\bH}(p_{B},\textbf{i}\mathbb{R}^{+})\geq \frac{\lsys}{4}$, by the triangle inequality we have that for any $z=(r,\theta)\in B(p_B; 1)$,
\beqar
\dist_{\bH}(z,\textbf{i}\mathbb{R^+})&\geq& \dist_{\bH}(p_B,\textbf{i}\mathbb{R^+})-\dist_{\bH}(p_B,z)\\
&\geq & \frac{\lsys}{4}-1\\
&\geq &\frac{\lsys}{8}.
\eeqar
Then by \eqref{i-dis} we know that
\[\sin(\theta)\leq 2e^{-\frac{\lsys}{8}}\]
which completes the proof.
\ep

Now we follow the same argument of the proof of \cite[Proposition 4.4]{W-inradius} to prove the following property with an effective leading constant.
\bpro \label{up-eff}
Let $X\in \sM_g$ with $\lsys\geq 8$. Then for any curve $\alpha \subset X$ with $\ell_{\alpha}(X)=\lsys$ there exists a uniform constant $C>0$ independent of $g$ such that
\[||\grad \ell_{\alpha}(X)||_{wp}^2\leq \frac{2}{\pi}\lsys\left(1+Ce^{-\frac{\lsys}{8}}  \right).\]
That is 
\[ ||\grad \sqrt{\ell_{\alpha}(X)}||_{wp}\leq \frac{1}{\sqrt{2\pi}}\sqrt{\left(1+Ce^{-\frac{\lsys}{8}}  \right)}.  \]
\epro

\begin{proof}
We will apply (\ref{Rie-f}) of Riera to finish the proof. First we know that
\[\lim_{u\to \infty} \frac{u \ln{\frac{u+1}{u-1}}-2}{u^{-2}}=\frac{2}{3}.\]

\noindent Similar as \cite[Equation $(4.5)$]{W-inradius}, since $\lsys\geq 8$, the quantity $u$ in Equation (\ref{Rie-f})  satisfies 
$u\geq \cosh(2)>1$. Thus, it follows by (\ref{Rie-f}) that there exists a uniform constant $C_2>0$ independent of $g$ such that 
\begin{eqnarray}
\langle \grad \ell_{\alpha},\grad \ell_{\alpha}\rangle_{wp}(X)\leq \frac{2}{\pi}\left(\ell_{\alpha}+C_{2} \displaystyle \sum_{B\in \left\{\left<A\right>\backslash\Gamma/\left<A\right>-id \right\}} e^{-2\dist_{\mathbb{H}}(\textbf{i}\mathbb{R}^+,B\circ (\textbf{i}\mathbb{R}^+))}\right). \nonumber
\end{eqnarray}

\noindent As introduced above one may choose $p_B=(r_B,\theta_B)\in B\circ (\textbf{i}\mathbb{R}^+)$ such that 
\ben
\item $\dist_{\bH}(p_{B},\textbf{i}\mathbb{R}^{+})\geq \frac{\lsys}{4}\geq 2$;

\item $1\leq r_{B}\leq e^{\lsys}$;

\item $\dist_{\mathbb{H}}(p_B, \textbf{i}\mathbb{R}^{+})=\dist_{\mathbb{H}}(B\circ (\textbf{i}\mathbb{R}^{+}),\textbf{i}\mathbb{R}^{+}).$
\een

\noindent Then, we have
\begin{eqnarray}\label{u-1}
\quad \quad ||\grad \ell_{\alpha}(X)||_{wp}^2 \leq  \frac{2}{\pi}\left(\ell_{\alpha}+C_{2} \displaystyle \sum_{B\in \left\{\left<A\right>\backslash\Gamma/\left<A\right>-id \right\}} e^{-2\dist_{\mathbb{H}}(\textbf{i}\mathbb{R}^+,p_B)}\right).
\end{eqnarray}

\noindent It is known from \cite[Lemma 2.4]{Wolpert08} or \cite[Lemma 2.1]{W-inradius} that the function  $e^{-2\dist_{\mathbb{H}}(\textbf{i}\mathbb{R}^+,z)}$ has the mean value property. More precisely, it follows by \cite[Lemma 2.4]{Wolpert08} or \cite[Lemma 2.1]{W-inradius} that there exists a uniform constant $C_3>0$ such that
\begin{eqnarray*}
 e^{-2\dist_{\mathbb{H}}(\textbf{i}\mathbb{R}^+,p_B)}\leq  C_3 \int_{B_{\mathbb{H}}(p_B;1)}e^{-2\dist_{\mathbb{H}}(z, \textbf{i}\mathbb{R^+})}\dArea(z).
\end{eqnarray*}

\noindent By Lemma \ref{unit-disjoint} we know that the geodesic balls $\{B_{\mathbb{H}}(p_B;1)\}_{B\in \left\{\left<A\right>\backslash\Gamma/\left<A\right>-id\right\}}$ are pairwise disjoint. Thus, we have
\begin{eqnarray}\label{u-2}
&& \sum_{B\in \left\{\left<A\right>\backslash\Gamma/\left<A\right>-id\right\}}  e^{-2\dist_{\mathbb{H}}(\textbf{i}\mathbb{R}^+,p_B)}\\
&\leq&  C_3 \sum_{B\in \left\{\left<A\right>\backslash\Gamma/\left<A\right>-id\right\}} \int_{B_{\mathbb{H}}(p_B;1)}e^{-2\dist_{\mathbb{H}}(z, \textbf{i}\mathbb{R^+})}\dArea(z) \nonumber\\ 
&=&C_3\displaystyle\int_{\bigcup_{ B\in \left\{\left<A\right>\backslash\Gamma/\left<A\right>-id\right\}}B_{\mathbb{H}}(p_B;1)}e^{-2\dist_{\mathbb{H}}(z, \textbf{i}\mathbb{R^+})}\dArea(z).\nonumber
\end{eqnarray}

\noindent It follows by \eqref{i-exp} and Lemma \ref{range-union} that
\bear \label{u-3}
&&\displaystyle\int_{\bigcup_{ B\in \left\{\left<A\right>\backslash\Gamma/\left<A\right>-id\right\}}B_{\mathbb{H}}(p_B;1)}e^{-2\dist_{\mathbb{H}}(z, \textbf{i}\mathbb{R^+})}\dArea(z) \\
&&\leq \int_{\sin (\theta)\leq 2 e^{-\frac{\lsys}{8}}} \int_{e^{-1}}^{e^{\lsys+1}}\sin^2{\theta}\dArea(z) \nonumber\\
&&= 2 \int_{0}^{\arcsin ( 2 e^{-\frac{\lsys}{8}})} \int_{e^{-1}}^{e^{\lsys+1}}\frac{\sin^2{\theta}}{r^2\sin^2{\theta}}rdrd\theta \nonumber\nonumber\\
&&=2\arcsin ( 2 e^{-\frac{\lsys}{8}})(\lsys+2) \nonumber\\
&&\leq C_4  \lsys e^{-\frac{\lsys}{8}} \nonumber
\eear
where $C_4>0$ is a uniform constant.

Thus, it follows by \eqref{u-1}, \eqref{u-2} and \eqref{u-3} that
\bear
||\grad \ell_{\alpha}(X)||_{wp}^2 \leq \frac{2}{\pi} \lsys \left(1+C_2C_3C_4 e^{-\frac{\lsys}{8}}\right).
\eear

Then the conclusion follows by choosing $$C=C_2C_3C_4>0$$ which is a uniform constant independent of $g$.
\end{proof}

\
\newline
Recall that
\[\sys(g)=\max_{X\in \sM_g}\lsys.\]
As $g\to \infty$, by Buser-Sarnak \cite{BS94} we know that $\sys(g)$ is uniformly comparable to $\ln(g)$. In particular $$\sys(g)\to \infty \ \emph{as} \ g\to \infty.$$

\noindent For any multicurve $\gamma \subset S_g$, we denote by $\sM_g^{\gamma}$ the stratum of $\sM_g$ whose pinching curves are $\gamma$. Wolpert in \cite{Wolpert08} applied Riera's formula \cite[Theorem 2]{Rie05} to give an upper bound for the \WP \ distance from any $X\in \sM_g$ to $\sM_g^{\gamma}$ in terms of the length $\ell_{\gamma}(X)$ of $\gamma$ at $X$. More precisely,
\bt[Wolpert, {\cite[Section 4]{Wolpert08}}]\label{w-d-up}
For any $X\in \sM_g$,
\[\dist_{wp}(X,\sM_g^{\gamma})\leq \sqrt{2\pi \ell_X(\gamma)}.\]
\et

Now we are ready to prove Theorem \ref{mt-3}.
\bt[=Theorem \ref{mt-3}]
The following limit holds:
\[\lim_{g\to \infty}\frac{\Inrad(\sM_g)}{\sqrt{\sys(g)}}=\sqrt{2\pi}\sim 2.5066.\]
\et

\bp[Proof of Theorem \ref{mt-3}]

For the upper bound, we follow the same argument as the proof of the upper bound of \cite[Theorem 1.1]{W-inradius}. For any hyperbolic surface $X \in \sM_g$, let $\alpha \subset X$ with $\ell_{sys}(X)=\ell_{\alpha}(X)$.
Let $\sM_g^{\alpha}$ be a stratum of $\overline{\sM_g}$ whose pinching curve is $\alpha$. Then it follows by Theorem \ref{w-d-up} that
\[\dist_{wp}(X, \sM_g^{\alpha})\leq \sqrt{2\pi \ell_{\alpha}(X)}\leq\sqrt{2\pi\sys(g)}\]
which implies that
\[\dist_{wp}(X,\p \sM_g)\leq \sqrt{2\pi\sys(g)}.\]
Since $X\in \sM_g$ is arbitrary, we have
\bear \label{up-2pi}
\limsup_{g\to \infty}\frac{\Inrad(\sM_g)}{\sqrt{\sys(g)}} \leq \sqrt{2\pi}\sim 2.5066.
\eear

For the lower bound, we follow a similar idea as in \cite{BB-20}. Let $X\in \sM_g$ with $$\lsys=\sys(g).$$ Let $c:[0, \Inrad(\sM_g)] \to \overline{\sM_g}$ be a \wep geodesic of arc-length parameter realizing $\Inrad(\sM_g)$, i.e., 
\ben
\item $c(0)=X$;

\item $c(t)\in \sM_g$ for all $t\in [0, \Inrad(\sM_g))$;

\item $c(\Inrad(\sM_g))\in \p \sM_g$.
\een  

\noindent For any fixed number $$T>8,$$ by continuity one may assume that $t_g \in (0, \Inrad(\sM_g))$ such that for large enough $g$,
\bear\label{low-T}
\min_{t\in [0,t_g]}\ell_{sys}(c(t))=\ell_{sys}(c(t_g))= T>8.
\eear 
By \cite[Lemma 3.4]{W-inradius} we know that $\ell_{sys}(c(\cdot))$ is piecewise smooth. So one may apply the Fundamental Theorem of Calculus and the Cauchy-Schwarz inequality to get
\bear \label{6-11}
\quad \left| \sqrt{\sys(g)}- \sqrt{T}\right|&=&\left| \sqrt{\lsys}- \sqrt{\ell_{sys}(c(t_g))}\right|\\
&=&\left| \int_{0}^{t_g} \langle \grad \sqrt{\ell_{sys}(c(t))}, c'(t) \rangle_{wp} dt \right| \nonumber\\
&\leq&  \int_{0}^{t_g} \left|\left| \grad \sqrt{\ell_{sys}(c(t))} \right|\right|_{wp} dt \nonumber \\
&\leq &  \frac{t_g}{\sqrt{2\pi}}\sqrt{\left(1+Ce^{-\frac{T}{8}}\right)} \nonumber 
\eear
where we apply Proposition \ref{up-eff} and \eqref{low-T} in the last inequality above, and the uniform constant $C>0$ is from Proposition \ref{up-eff}. Clearly we have $$t_g \leq \Inrad(\sM_g).$$ Recall that $\sys(g) \to \infty$ as $g\to \infty$. Thus, it follows by \eqref{6-11} that
\bear
\liminf_{g\to \infty}\frac{\Inrad(\sM_g)}{\sqrt{\sys(g)}}&\geq& \liminf_{g\to \infty}\frac{t_g}{\sqrt{\sys(g)}} \nonumber\\
&\geq& \frac{\sqrt{2\pi}}{\sqrt{\left(1+Ce^{-\frac{T}{8}}\right)}}.\nonumber
\eear
Since $T\geq 8$ is arbitrary, let $T\to \infty$ we get
\bear\label{low-2pi}
\liminf_{g\to \infty}\frac{\Inrad(\sM_g)}{\sqrt{\sys(g)}}\geq \sqrt{2\pi}\sim 2.5066.
\eear

Then the conclusion follows by \eqref{up-2pi} and \eqref{low-2pi}.
\ep

\br
The argument for Theorem \ref{mt-3} also works for the moduli space $\sM_{g,n}$ of \RS s with punctures where $n>0$. We only consider the closed case in this paper for simplicity.
\er

\br
The argument of Theorem \ref{mt-3} above highly relies on large genus. Bromberg and Bridgeman in \cite{BB-20} proved the following surprising result: for all $g,n$ with $3g+n-3<0$, $$\Inrad(\sM_{g,n})\geq 0.94 \sqrt{2\pi\sys(g,n)}$$
where $\sys(g,n)=\max_{X\in \sM_{g,n}}\lsys$.  
\er

Similarly, we define the \wep \emph{inradius} $\Inrad(\sT_g)$ of the Teichm\"uller space $\sT_g$ as
\[\Inrad(\sT_g):=\max_{X \in \sT_g}\dist_{wp}(X, \p \sT_g)\]
where $\p \sT_g$ is the boundary of $\sT_g$ consisting of nodal surfaces. The proof of Theorem \ref{mt-3} also gives that
\begin{corollary} The following limit holds:
\[\lim_{g\to \infty}\frac{\Inrad(\sT_g)}{\sqrt{\sys(g)}}=\sqrt{2\pi}.\]
\end{corollary}

\section{Proof of Theorem \ref{mt-4}}\label{sec-mt4}

For any $L>0$ which may depend on the genus $g$, recall that for any $X\in \mathcal{MC}(\leq L)$ there exists a pants decomposition $\mathcal{P}$ of $X$ such that the length satisfies $$\max_{\alpha \in \mathcal{P}}\ell(\alpha)\leq L.$$ The following lemma tells that the largest radius of embedded hyperbolic disk in $X$ is bounded above by a function of $L$. More precisely,
\bl \label{mcl-up}
For any $X\in \mathcal{MC}(\leq L)$,
\[\max_{p \in X} \inj_X(p)<\frac{L}{2}+\ln(6).\]
\el 
\bp
For any $p \in X \in  \mathcal{MC}(\leq L)$, one may assume that $p \in P$ where $P$ is a pant whose three boundary closed geodesics all have length $\leq L$. Take two copies of $P$ and we double them into a closed hyperbolic surface $X_2$ of genus $2$ (here one may take any twist along these three closed geodesics). In $X_2$, by \eqref{up-inj} we have 
\[\inj_{X_2}(p)\leq \ln(6).\] 
Let $B_{X_2}(p;\inj_{X_2}(p)):=\left\{q\in X_2; \ \dist(q, p)<\inj_{X_2}(p)\right\}$ be the hyperbolic open disk in $X_2$ centered at $p$ of radius $\inj_{X_2}(p)$. We finish the proof by considering the following two cases.

\emph{Case $(1)$. $B_{X_2}(p;\inj_{X_2}(p)) \cap \p P =\emptyset$.} For this case we clearly have
\bear
\inj_X(p)=\inj_{X_2}(p)\leq \ln(6).
\eear

\emph{Case $(2)$. $B_{X_2}(p;\inj_{X_2}(p)) \cap \p P \neq \emptyset$.} For this case, one may assume $$p_0 \in B_{X_2}(p;\inj_{X_2}(p))\cap \alpha$$ for some component $\alpha$ of $\p P$. We parametrize $\alpha$ such that $\alpha(0)=\alpha(\ell(\alpha))=p_0$. Let $\sigma:[0,\dist(p,p_0)]\to P\subset X$ be a shortest geodesic of $X$ joining $p$ and $p_0$. Consider the closed curve $\sigma'$ based at $p$ as following.
{\small{
\beqar
\sigma'(s):=
\begin{cases}
\sigma(s), \quad \emph{$0\leq s\leq \dist(p,p_0)$}\\
\alpha(s-\dist(p,p_0)), \quad \emph{$\dist(p,p_0)\leq s \leq \dist(p,p_0)+\ell(\alpha)$}\\
\sigma(2\dist(p,p_0)+\ell(\alpha)-s), \quad \emph{$\dist(p,p_0)+\ell(\alpha) \leq s \leq 2\dist(p,p_0)+\ell(\alpha)$}
\end{cases}
\eeqar  
}}
This closed curve $\sigma' \subset P\subset X$ is freely homotopic to $\alpha$. So it is nontrivial. Thus, we have
\bear
\inj_X(p)<\frac{\ell(\sigma')}{2}=\dist(p,p_0)+\frac{\ell(\alpha)}{2}\leq \ln(6)+\frac{L}{2}.
\eear

Then the conclusion follows by these two cases.
\ep

Now we recall the following result of Mirzakhani which roughly says that as $g\to \infty$, almost all hyperbolic surfaces of genus $g$ contain an embedded hyperbolic disk of radius $\frac{\ln(g)}{6}$. More precisely,
\bt[Mirzakhani, {\cite[Theorem 4.5]{Mirz13}}] \label{MM13}
\[\lim_{g\to \infty}\frac{\Vol\left(\left\{X\in \sM_g;\ \max_{p \in X} \inj_X(p)\geq \frac{\ln(g)}{6}\right\}\right)}{\Vol(\sM_g)}=1.\]
\et

Now we are ready to prove Theorem \ref{mt-4}.
\bt[= Theorem \ref{mt-4}]
For any small enough $\epsilon>0$ and if $L_g=\epsilon \ln(g)$, then 
\[\lim_{g\to \infty}\frac{\Vol\left(\left\{X\in \sM_g; \ \dist_{wp}\left(X,\mathcal{MC}(\leq L_g)\right)> 1.0511\sqrt{\ln(g)}\right\}\right)}{\Vol(\sM_g)}=1\]
where $\Vol(\cdot)$ is the \WP \ volume.
\et

\bp [Proof of Theorem \ref{mt-4}]
Set $$\mathcal{MR}:=\left\{X\in \sM_g; \ \max_{p \in X} \inj_X(p)\geq \frac{\ln(g)}{6} \right\}.$$
For any $X\in \mathcal{MR}$ and $Y\in \mathcal{MC}(\leq L_g)$, one may let $p \in S_g$ such that 
\[\inj_X(p)\geq \frac{\ln(g)}{6}.\]

Then it follows by Theorem \ref{mt} and Lemma \ref{mcl-up} that
\beqar
\dist_{wp}(X,Y)&\geq& \frac{\sqrt{\inj_X(p)}-\sqrt{\inj_Y(p)}}{0.3884}\\
&\geq & \frac{\sqrt{\ln(g)}-\sqrt{3L_g+6\ln(6)}}{0.3884\sqrt{6}}.
\eeqar
Since $\frac{1}{0.3884\sqrt{6}}\sim1.051102>1.0511$ and $L_g=\epsilon \ln(g)$, we have that for large enough $g$ and small enough $\epsilon>0$,
\bear
\dist_{wp}(X,Y)>1.0511\sqrt{\ln(g)}.
\eear
Which implies that for large enough $g$,
\[\mathcal{MR}\subset \left\{X\in \sM_g; \ \dist_{wp}\left(X,\mathcal{MC}(\leq L_g)\right)> 1.0511\sqrt{\ln(g)}\right\}.\]
Then the conclusion follows by Theorem \ref{MM13}.
\ep

\br
For $L>0$, we say $X\in \sM_g$ has \emph{total pants length at least $L$} if for any pants decomposition $\mathcal{P}$ of $X$, $\sum_{\alpha\in \mathcal{P}}\ell(\alpha)\geq L$. Guth, Parlier and Young in \cite{GPY11} showed that for any $\eps>0$, 
\[\lim_{g\to \infty}\frac{\Vol\left(\left\{X\in \sM_g;\ \emph{$X$ has total pants length at least $g^{\frac{7}{6}-\eps}$}\right\}\right)}{\Vol(\sM_g)}=1.\]
Clearly we have that for any $X\in \mathcal{MC}(\leq L_g)$, the hyperbolic surface $X$ has total pants length at most $(3g-3)L_g$. We do \emph{not} know too much information on the least total pants length of $Y\in \sM_g$ with $\dist_{wp}(Y, \mathcal{MC}(\leq L_g))> 1.0511\sqrt{\ln(g)}$.
\er

\bibliographystyle{amsalpha}
\bibliography{wp}

\providecommand{\bysame}{\leavevmode\hbox to3em{\hrulefill}\thinspace}
\providecommand{\MR}{\relax\ifhmode\unskip\space\fi MR }
\providecommand{\MRhref}[2]{%
  \href{http://www.ams.org/mathscinet-getitem?mr=#1}{#2}
}
\providecommand{\href}[2]{#2}
\begin{thebibliography}{GPY11}

\bibitem[BB22]{BB-20}
Martin {Bridgeman} and Kenneth {Bromberg}, \emph{{Strata Separation for the
  Weil-Petersson Completion and Gradient Estimates for Length Functions}},
  Journal of Topology and Analysis, {doi: 10.1142/S1793525321500667,} (2022).

\bibitem[BBB19]{BBB19}
Martin Bridgeman, Jeffrey Brock, and Kenneth Bromberg, \emph{Schwarzian
  derivatives, projective structures, and the {W}eil-{P}etersson gradient flow
  for renormalized volume}, Duke Math. J. \textbf{168} (2019), no.~5, 867--896.

\bibitem[Bro03]{Brock03}
Jeffrey~F. Brock, \emph{The {W}eil-{P}etersson metric and volumes of
  3-dimensional hyperbolic convex cores}, J. Amer. Math. Soc. \textbf{16}
  (2003), no.~3, 495--535 (electronic).

\bibitem[BS94]{BS94}
P.~Buser and P.~Sarnak, \emph{On the period matrix of a {R}iemann surface of
  large genus}, Invent. Math. \textbf{117} (1994), no.~1, 27--56, With an
  appendix by J. H. Conway and N. J. A. Sloane.

\bibitem[Bus10]{Buser10}
Peter Buser, \emph{Geometry and spectra of compact {R}iemann surfaces}, Modern
  Birkh\"auser Classics, Birkh\"auser Boston, Inc., Boston, MA, 2010, Reprint
  of the 1992 edition.

\bibitem[BW21]{BW-curv}
Martin Bridgeman and Yunhui Wu, \emph{Uniform bounds on harmonic {B}eltrami
  differentials and {W}eil-{P}etersson curvatures}, J. Reine Angew. Math.
  \textbf{770} (2021), 159--181.

\bibitem[CP12]{CP12}
William Cavendish and Hugo Parlier, \emph{Growth of the {W}eil-{P}etersson
  diameter of moduli space}, Duke Math. J. \textbf{161} (2012), no.~1,
  139--171.

\bibitem[GPY11]{GPY11}
Larry Guth, Hugo Parlier, and Robert Young, \emph{Pants decompositions of
  random surfaces}, Geom. Funct. Anal. \textbf{21} (2011), no.~5, 1069--1090.

\bibitem[IT92]{IT92}
Y.~Imayoshi and M.~Taniguchi, \emph{An introduction to {T}eichm\"uller spaces},
  Springer-Verlag, Tokyo, 1992, Translated and revised from the Japanese by the
  authors.

\bibitem[KM18]{KM18}
Sadayoshi Kojima and Greg McShane, \emph{Normalized entropy versus volume for
  pseudo-{A}nosovs}, Geom. Topol. \textbf{22} (2018), no.~4, 2403--2426.

\bibitem[Mir13]{Mirz13}
Maryam Mirzakhani, \emph{Growth of {W}eil-{P}etersson volumes and random
  hyperbolic surfaces of large genus}, J. Differential Geom. \textbf{94}
  (2013), no.~2, 267--300.

\bibitem[Rie05]{Rie05}
Gonzalo Riera, \emph{A formula for the {W}eil-{P}etersson product of quadratic
  differentials}, J. Anal. Math. \textbf{95} (2005), 105--120.

\bibitem[RT18a]{RT18}
Melanie Rupflin and Peter~M. Topping, \emph{Horizontal curves of hyperbolic
  metrics}, Calc. Var. Partial Differential Equations \textbf{57} (2018),
  no.~4, Paper No. 106, 17.

\bibitem[RT18b]{RT18-JDG}
\bysame, \emph{Teichm\"{u}ller harmonic map flow into nonpositively curved
  targets}, J. Differential Geom. \textbf{108} (2018), no.~1, 135--184.

\bibitem[Sch13]{Schl13}
Jean-Marc Schlenker, \emph{The renormalized volume and the volume of the convex
  core of quasifuchsian manifolds}, Math. Res. Lett. \textbf{20} (2013), no.~4,
  773--786.

\bibitem[Sch20]{Schl19}
\bysame, \emph{{Volumes of quasifuchsian manifolds}}, Surveys in Differential
  Geometry (2020), \emph{to appear}.

\bibitem[Teo09]{Teo09}
Lee-Peng Teo, \emph{The {W}eil-{P}etersson geometry of the moduli space of
  {R}iemann surfaces}, Proc. Amer. Math. Soc. \textbf{137} (2009), no.~2,
  541--552.

\bibitem[Tro92]{Tromba-book}
Anthony~J. Tromba, \emph{Teichm\"{u}ller theory in {R}iemannian geometry},
  Lectures in Mathematics ETH Z\"{u}rich, Birkh\"{a}user Verlag, Basel, 1992,
  Lecture notes prepared by Jochen Denzler.

\bibitem[Wol87]{Wolpert87}
Scott~A. Wolpert, \emph{Geodesic length functions and the {N}ielsen problem},
  J. Differential Geom. \textbf{25} (1987), no.~2, 275--296.

\bibitem[Wol08]{Wolpert08}
\bysame, \emph{Behavior of geodesic-length functions on {T}eichm\"uller space},
  J. Differential Geom. \textbf{79} (2008), no.~2, 277--334.

\bibitem[Wol10]{Wolpert-book}
\bysame, \emph{Families of {R}iemann surfaces and {W}eil-{P}etersson geometry},
  CBMS Regional Conference Series in Mathematics, vol. 113, Published for the
  Conference Board of the Mathematical Sciences, Washington, DC; by the
  American Mathematical Society, Providence, RI, 2010.

\bibitem[Wu19]{W-inradius}
Yunhui Wu, \emph{Growth of the {W}eil-{P}etersson inradius of moduli space},
  Annales de l'Institut Fourier \textbf{69} (2019), no.~3, 1309--1346 (en).

\bibitem[Wu20]{W-sys}
\bysame, \emph{{Systole functions and {W}eil-{P}etersson geometry}},
  \emph{preprint}.

\bibitem[WW18]{Wolf-W-1}
Michael Wolf and Yunhui Wu, \emph{Uniform bounds for {W}eil-{P}etersson
  curvatures}, Proc. Lond. Math. Soc. (3) \textbf{117} (2018), no.~5,
  1041--1076.

\end{thebibliography}

\end{document}